\documentclass[a4paper,regno,11pt]{amsart}
\usepackage{latexsym, enumerate}

\usepackage[pdftex]{graphicx}

\usepackage{xcolor}
\usepackage[percent]{overpic}

\usepackage{amsmath, amsthm}
\usepackage{amsfonts}
\usepackage{amssymb}
\usepackage{fancyhdr}
\usepackage{setspace}
\usepackage{array}
\usepackage{color}
\usepackage{xmpmulti}
\usepackage{enumitem}

\newcommand{\ds}{\displaystyle }
\newtheorem{thm}{Theorem}[section]
\newtheorem{rem}[thm]{Remark}

\begin{document}

\title{Macroscopic modeling and simulations of room evacuation}

\author[M. Twarogowska, P. Goatin, R. Duvigneau]{M. Twarogowska$^1$ \and P. Goatin$^2$ \and R. Duvigneau$^2$}

\thanks{$^1$ Istituto per le Applicazioni del Calcolo ``Mauro Picone'', Consiglio Nazionale delle Ricerche,  via dei Taurini 19, I-00185 Roma, Italy
({\tt mtwarogowska@gmail.com}).}

\thanks{$^2$ INRIA Sophia Antipolis - M\'editerran\'ee, 
OPALE Project-Team,
2004, route des Lucioles -- BP 93, 
06902 Sophia Antipolis Cedex, France
({\tt paola.goatin@inria.fr}, {\tt regis.duvigneau@inria.fr})
}         
         
\begin{abstract}
We analyze numerically two macroscopic models of crowd dynamics: the classical Hughes model and the second order model being an extension to pedestrian motion of the Payne-Whitham vehicular traffic model. The desired direction of motion is determined by solving an eikonal equation with density dependent running cost, which results in minimization of the travel time and avoidance of congested areas. We apply a mixed finite volume-finite element method to solve the problems and present error analysis for the eikonal solver, gradient computation and the second order model yielding a first order convergence. We show that Hughes' model is incapable of reproducing complex crowd dynamics such as stop-and-go waves and clogging at bottlenecks. Finally, using the second order model, we study numerically the evacuation of pedestrians from a room through a narrow exit.
\end{abstract}

\thanks{ {\emph{keywords and phrases: }} macroscopic models,  crowd dynamics, evacuation,  Braess paradox }

\maketitle

\section{Introduction}
Crowd dynamics has recently attracted the interests of a rapidly increasing number of scientists. Analytical and numerical analysis are effective tools to investigate, predict and simulate complex behaviour of pedestrians, and numerous engineering applications welcome the support of mathematical modelling. Growing population densities combined with easier transport lead to greater accumulation of people and increase risk of life threatening situations. Transport systems, sports events, holy sites or fire escapes are just few examples where uncontrolled behaviour of a crowd may end up in serious injuries and fatalities. In this field, pedestrian traffic management is aimed at designing walking facilities which follow optimal requirements regarding flow efficiency, pedestrians comfort and, above all, security and safety. 

From a mathematical point of view, a description of human crowds is strongly non standard due to the intelligence and decision making abilities of pedestrians. Their behaviour depends on the physical form of individuals and on the purpose and conditions of their motion. In particular, pedestrians walk with the most comfortable speed \cite{Buchmueller}, tend to maintain a preferential direction towards their destination and avoid congested areas. On the contrary, in life threatening circumstances, nervousness make them move faster \cite{predtechenskii1978planning}, push others and follow the crowd instead of looking for the optimal route \cite{Keating}. As a consequence, critical crowd conditions appear such as ''freezing by heating'' and ''faster is slower'' phenomena \cite{Helbing_freezing, Helbing2002}, stop-and-go waves, transition to irregular flow \cite{Helbing_Johansson_Abideen} and arching and clogging at bottlenecks  \cite{predtechenskii1978planning}. 

In order to describe this complex crowd dynamics, numerous mathematical models have been introduced, belonging to two fundamentally distinct approaches: microscopic and macroscopic. In the microscopic framework pedestrians are treated as individual entities whose trajectories are determined by physical and social laws. Examples of microscopic models are the social force model \cite{Helbing95},  cellular automata models \cite{Muramatsu, Dijkstra}, AI-based models \cite{Gopal}. Macroscopic description treats the crowd as a continuum medium characterized by averaged quantities such as density and mean velocity. The first modelling attempt is due to Hughes \cite{Hughes2002} who defined the crowd as a ''thinking fluid'' and described the time evolution of its density using a scalar conservation law. Current macroscopic models use gas dynamics equations \cite{Bellomo, Jiang2010}, gradient flow methods \cite{Maury}, non linear conservation laws with non classical shocks \cite{ColomboRosini} and time evolving measures \cite{Tosin}. 
At an intermediate level, kinetic models derive evolution equations for the probability distribution functions of macroscopic variables directly from microscopic
interaction laws between individuals, see for example \cite{Bellomo_Bellouquid} and \cite{Degond_etal} and references therein. Also, recently introduced approaches include micro-macro coupling of time-evolving measures \cite{CristianiPiccoliTosin} and mean-field games \cite{LachapelleWolfram}. These models are good candidates to capture the effects of individual behavior on the whole system.\newline

In this paper we shall analyze and compare two macroscopic models describing the time evolution of the density of pedestrians. The first one, introduced by Hughes \cite{Hughes2002}, consists of a mass conservation equation supplemented with a phenomenological relation between the speed and the density of pedestrians. The second one involves mass and momentum balance equations so is of second order type. It was proposed by Payne and Whitham \cite{Payne1971, Whitham1974} to describe vehicular traffic and adopted to describe pedestrian motion by Jiang et al. \cite{Jiang2010}. It consists of the two-dimensional Euler equations with a relaxation source term. In both models, the pedestrians' optimal path is computed using the eikonal equation as was proposed by Hughes \cite{Hughes2002}. 

In order to simulate realistic behaviour we consider two dimensional, continuous walking domains with impenetrable walls and exits as pedestrians' destination. To our knowledge the only available results using Hughes' model concern simulations of flow of pedestrians on a large platform with an obstacle in its interior \cite{Huang, Jiang_Liu}. In the case of the second order model Jiang et al. \cite{Jiang2010} considered the same setting and showed numerically the formation of stop-and-go waves. However, none of the above works analyzed complex crowd dynamics. Behaviour at bottlenecks and evacuation process was not considered in any of the previous works.  \newline

The first aim of this paper is to provide a more detailed insight into the properties of macroscopic models of pedestrian motion. In particular, we compare Hughes' model and the second order model analyzing the formation of stop-and-go waves and flows through bottlenecks. Our simulations suggest that Hughes' model is incapable of reproducing neither such waves nor clogging at a narrow exit. It appears to be also insensitive to the presence of obstacles placed in the interior of the walking domain, which can be crucial in the study of evacuation.  This is why in the second part of the paper we restrict ourselves only to the second order model 

We focus on the study of the evacuation of pedestrians through a narrow exit. This problem is an important safety issue because of arching and clogging appearing in front of the exit, which can interrupt the outflow and result in crushing of people under the pressure or the crowd. Experimental studies are rare due to the difficulties in reproducing realistic panic behaviour \cite{Kretz, Seyfried_bottleneck, HoogendoornDaamen}, while numerical simulations are available mainly in the microscopic framework. For example Helbing et al. \cite{Helbing_Nature, Helbing_Johansson2009} analyzed the evacuation of two hundred people from a room through a narrow door and in \cite{Helbing_Buzna2005} the issue of optimal design of walking facilities was addressed with genetic algorithms. 
At first we show the dependence of the solutions on different parameters of the model. More precisely, we consider the effect on the evacuation of the strength of the interpersonal repealing forces and the desired speed of pedestrians. Both of these parameters may indicate the nervousness and the level of panic of pedestrians. 

In order to improve evacuation, Hughes \cite{Hughes2003} suggested that suitably placed obstacles can increase the flow through an exit. This idea is an inversion of the Braess paradox \cite{Braess1968, Braess2005}, which was formulated for traffic flows and states that adding extra capacity to a network can in some cases reduce the overall performance. In the case of crowd dynamics, placing an obstacle may be seen intuitively as a worse condition. Nevertheless, it is expected to lower the internal pressure between pedestrians and their density in front of the exit and as a result preventing from clogging. This phenomenon has been studied experimentally in case of granular materials by Zuriguel et al. \cite{Zuriguel} who analyzed the outflow of grains from a silo and found out the optimal height above the outlet of an obstacle which reduces the blocking of the flow by a factor of one hundred. In case of pedestrians, to our knowledge, so far this problem has been studied only numerically. Helbing et al. \cite{Helbing2002} using the social force model observed that a single column placed in front of the exit decreases the pressure between the column and the door and may prevent from clogging. In the same framework, different shapes and placements of obstacles were studied in \cite{Escobar} with an indication of the formation of the so called ''waiting zone'' in front of the exit. Frank and Dorso in \cite{Frank_Dorso} studied the effects of a column and a longitudinal panel assuming in the social force model that pedestrians change their direction away from an obstacle until the exit becomes visible. 
  
Following the idea of Hughes \cite{Hughes2003}, we try to improve the evacuation of pedestrians using properly tuned obstacles placed in front of the exit. Motivated by the numerical simulations in which clogging appears when a large group of pedestrians reached the exit simultaneously, we give an example of a system of five circular columns arranged in the shape of a triangle opened towards the exit. We show that this system of obstacles effectively creates an area with lower density in front of the door and reduces the clogging.\newline

This paper is organized as follows: in Section~\ref{sec:models} we explain in detail macroscopic models and in Section~\ref{sec:numericalScheme} we describe numerical approximation of the models. Section~\ref{sec:results} is devoted to the numerical results. At first we present error analysis and comparison between the two macroscopic models. Then we analyze the evacuation of pedestrians from a room. 

\section{Macroscopic model of pedestrian flow}\label{sec:models}         

\subsection{Equations}

We consider a two dimensional connected domain $\Omega\subset\mathbb{R}^2$ corresponding to some walking facility. It is equipped with an exit which models the destination of the crowd motion and can contain obstacles. The boundary of the domain $\Omega$ is composed of the outflow boundary $\Gamma_{o}$ and the wall $\Gamma_{w}$, which, as obstacles, is impenetrable for the pedestrians. In this setting we consider a macroscopic model introduced by Payne-Whitham for vehicular traffic flow in \cite{Payne1971, Whitham1974} and by Jiang et al. in \cite{Jiang2010} to describe crowd dynamics. The model derives from fluid dynamics and consists of mass and momentum balance equations with source term. Denoting by $\rho$ the density of pedestrians and by $\vec{v}$ their mean velocity the model reads

      \begin{equation}\label{eq:mainSystem}
            \left\{
                 \begin{array}{l}
                   \rho_{t}+\textrm{div}(\rho\vec{v})=0,\\
                   (\rho\vec{v})_{t}+\textrm{div}(\rho\vec{v}\otimes\vec{v})=\vec{\mathcal{A}}(\rho,\vec{v}),
                 \end{array}\right.
      \end{equation}
where $\vec{\mathcal{A}}{(\rho,\vec{v})}$ describes the average acceleration caused by internal driving forces and motivations of pedestrians. More precisely, it consists of a relaxation term towards a desired velocity and the internal pressure preventing from overcrowding
      \begin{equation}\label{eq:vectorF}
            \vec{\mathcal{A}}{(\rho,\vec{v})}= \frac{1}{\tau}\left(\rho V(\rho)\vec{\mu}-\rho\vec{v}\right)-\nabla P(\rho).\color{red}{}
      \end{equation}
The unit vector $\vec{\mu}=\vec{\mu}(\rho(x,t))$ describes the preferred direction pointing the objective of the movement of pedestrians and will be defined in the next section. The function $V(\rho)$ characterizes how the speed of pedestrians changes with density. Various speed-density relations are available in the literature, see \cite{Buchmueller}. For our simulations we choose the exponential dependence
\begin{equation}\label{eq:speed}
      V(\rho) = v_{\max}{e}^{-\alpha\left(\frac{\rho}{\rho_{\max}}\right)^2},
\end{equation}
where $v_{\max}$ is a free flow speed, $\rho_{\max}$ is a congestion density at which the motion is hardly possible and $\alpha$ is a positive constant. The parameter $\tau$ in (\ref{eq:vectorF}) is a relaxation time describing how fast pedestrians correct their current velocity to the desired one. The second term in (\ref{eq:vectorF}) models a repulsive force modeling the volume filling effect and is given by the power law for isentropic gases
\begin{equation}\label{eq:pressure}
      P(\rho)=p_{0}\rho^\gamma,\qquad p_{0}>0,\quad\gamma>1. 
\end{equation}
\begin{rem}Model (\ref{eq:mainSystem}) is referred to as a second order model as it consists of mass and momentum balance equations completed with a phenomenological law describing the acceleration. A simpler system, which is a first order model, was introduced by Hughes \cite{Hughes2002, Hughes2003}. It is composed of a scalar conservation law
\begin{equation}\label{eq:Hughes}
      \rho_{t}+\textrm{div}\vec{F}(\rho)=0,
\end{equation}
where $\vec{F}(\rho)=\rho V(\rho)\vec{\mu}$, closed by a speed-density relation $V(\rho)$ given by \eqref{eq:speed}. 
\end{rem}

\subsection{Desired velocity}
The models (\ref{eq:mainSystem}), \eqref{eq:Hughes} have to be completed by defining the vector field $\vec{\mu}$. Following the works of Hughes,  we assume that the pedestrians movement is opposite to the gradient of a scalar potential $\phi$, that is
\begin{equation}\label{eq:directionVector}
      \vec{\mu}=-\frac{\nabla\phi}{||\nabla\phi||}.
\end{equation}
The potential $\phi$ corresponds to an instantaneous travel cost which pedestrians want to minimize and is determined by the eikonal equation 
\begin{equation}\label{eq:eikonal}
      \left\{
           \begin{array}{lcc}
              |\nabla\phi| = c(\rho)&\textrm{in}&\Omega \\
              \phi=0&\textrm{on}&\Gamma_{\textrm{o}}
           \end{array}\right.,
\end{equation}
 where  $c(\rho)$ is a density dependent cost function increasing with $\rho$. In the simplest case we could prescribe $c(\rho)=1$, which gives the potential $\phi(x)=\textrm{dist}(x,\Gamma_{o})$ in the case of convex domains. Pedestrians want to minimize the path towards their destination but temper the estimated travel time by avoiding high densities. The behaviour can be expressed by the ''density driven'' rearrangement of the equipotential curves of $\phi$ using the following cost function \cite{Hughes2002}
\begin{equation}\label{eq:runningCost}
      c(\rho)=\frac{1}{V(\rho)}.
\end{equation}

\begin{rem} Instead of coupling the mass and momentum balance laws with an eikonal equation, another possible approach has been recently introduced
in \cite{Maury, Maury_etal}: the transport equation is interpreted as a gradient flow in the Wasserstein space, which has the advantage of providing existence
results despite the non-smooth setting. 
\end{rem}

\section{Numerical scheme}\label{sec:numericalScheme}
Let us now present the numerical scheme on unstructured triangular mesh that we used to perform numerical simulations. The model of pedestrian flow couples equations of different nature, i.e. a two dimensional non-linear system of conservation laws with sources, coupled with the eikonal equation through the source term. In this section we describe a finite volume scheme built on dual cells for systems (\ref{eq:mainSystem}) and \eqref{eq:Hughes} and a finite element method based on the variational principle for problem (\ref{eq:eikonal}). The numerical simulations are carried out using the multidisciplinary platform NUM3SIS developed at Inria Sophia Antipolis \cite{num3sis, num3sisOpale}. 

\subsection{Finite volume schemes for the macroscopic models}
The models of pedestrian motion (\ref{eq:mainSystem}) and \eqref{eq:Hughes} can be put in the form
\begin{equation}\label{eq:mainSystemMatrix} 
      U_{t}+\textrm{div}\vec{F}(U)=S(U),
\end{equation}
where in the case of the second order model (\ref{eq:mainSystem}) $U=\left(\rho, \rho\vec{v}\right)^{T}$ denotes the unknowns vector, density and momentum, and 
\begin{displaymath}
\vec{F}(U)=
          \left(
               \begin{array}{c}
                  F^{\rho} \\
                  F^{\rho \vec{v}}
               \end{array}
          \right)
         =\left(
               \begin{array}{c}
                  \rho \vec{v}\\
                  \rho\vec{v}\otimes\vec{v}+P(\rho)
               \end{array}
          \right),
\quad S(U)=
          \left(
               \begin{array}{c}
                  0\\ 
                  \frac{1}{\tau}\left(\rho V(\rho)\vec{\mu}-\rho\vec{v}\right)
               \end{array}
          \right).
\end{displaymath}
For the first order model \eqref{eq:Hughes} we take $U=\rho$, $\vec{F}(U)=\rho V(\rho)\vec{\mu}$ and $S(U)=0$. 

According to the framework of finite volume schemes, we decompose the domain $\Omega$ into $N$ non overlapping, finite volume cells $C_{i}$, $i=1,...,N$, given by dual cells centered at vertices of the triangular mesh. For each cell $C_{i}$ we consider a set of $N_{i}$ neighbouring cells $C_{ij}$, $j=1,...,N_{i}$. By $e_{ij}$ we denote the face between $C_{i}$ and $C_{j}$, $|e_{ij}|$ its length and $\vec{n}_{ij}$ is a unit vector normal to the $e_{ij}$ pointing from the center of the cell $C_{i}$ towards the center of the cell $C_{j}$. The solution $U$ of the system (\ref{eq:mainSystemMatrix}) on a cell $C_{i}$ is approximated by the cell average of the solution at time $t>0$, that is
\begin{displaymath}
      \displaystyle{U_{i}=\frac{1}{|C_{i}|}\int_{C_{i}}U(x,t)dx.} 
\end{displaymath}
A general semi-discrete finite volume scheme for (\ref{eq:mainSystemMatrix}) can be defined as
\begin{equation}\label{eq:semiDiscreteScheme}
      \frac{d}{dt}U_{i}=-\frac{1}{|C_{i}|}\sum_{j=1}^{N_{i}}|e_{ij}|\mathcal{F}(U_i,U_j,\vec{n}_{ij})+S(U_{i}),
\end{equation}
where $\mathcal{F}(U_i,U_j,\vec{n}_{ij})$ is a numerical flux function. The spatial discretization of the source term $S(U_{i})$ is treated by a pointwise approximation $
S(U_{i})=\left( 0, \left(\rho_{i} V(\rho_{i})\vec{\mu}_{i}-\rho_{i}\vec{v}_{i}\right)/\tau\right)^T.$

In order to obtain a numerical approximation using a finite volume scheme \eqref{eq:semiDiscreteScheme} we have to compute numerical fluxes $\mathcal{F}(U_i,U_j,\vec{n}_{ij})$ across the face $e_{ij}$ between control cells $C_i$ and $C_j$ along the direction $\vec{n}_{ij}$.   Despite the fact that the model is two dimensional, these fluxes are computed using a one-dimensional approximation. 

The homogeneous part of the model \eqref{eq:mainSystem} coincides with the isentropic gas dynamics system for which many solvers are available, (see \cite{toro}). However, the occurrence of vacuum may cause instabilities and not all of them preserve non negativity of the density. We use the first order HLL approximate Riemann solver \cite{HLL}. It assumes that the solution consists of three constant states separated by two waves with speeds $\sigma_{L}$ and $\sigma_R$ corresponding respectively to the slowest and fastest signal speeds.  It is ''positivity preserving'' under certain conditions on the above numerical wave speeds \cite{Einfeldt} that is
\begin{equation}\label{eq:waveSpeeds}
      \begin{array}{l}
         \sigma_{L}=\min(v_{L}^{n}-s_{L},\bar{v}_{Roe}-\bar{s}),\\
         \sigma_{R}=\max(\bar{v}_{Roe}+\bar{s},v_{R}^{n}+s_{R}),
      \end{array}
\end{equation}
where $s=\sqrt{P'(\rho)}$ is the sound speed, $v^{n}$ is the normal component of the velocity and $\bar{v}_{Roe},\bar{s}$ are averaged Roe velocity and sound speed respectively. 

For the numerical function in the case of the first order model we use the Lax-Friedrichs flux
\begin{subequations}\label{eq:LaxFriedrichsFlux}
      \begin{equation}
            \mathcal{F}\left(\rho_{i},\rho_{j},\vec{n}_{ij}\right)=\frac{1}{2}\left[F(\rho_{i})\cdot\vec{n}_{ij}+F(\rho_{j})\cdot\vec{n}_{ij}-\xi(\rho_{j}-\rho_{i})\right],
      \end{equation}
with the numerical viscosity coefficient $\xi$ given by
      \begin{equation}
            \xi=\max_{l=i,j}\left|\frac{d}{d\rho}F(\rho_{l})\right|=\max_{l=i,j}\left|\frac{d}{d\rho}\left(\rho_{l} V(\rho_{l})\vec{\mu}_{l}(\rho)\right)\right|=\max_{l=i,j}\left|\frac{d}{d\rho}\left(\rho_{l} V(\rho_{l})\right)\right|.
      \end{equation}
\end{subequations}
The last equality is justified by the fact that $\vec{\mu}$ is a unit vector.

\subsection{Fully discrete schemes}
The difficulty in the time discretization of equation (\ref{eq:semiDiscreteScheme}) lies in the non linear coupling of the models with the eikonal equation (\ref{eq:eikonal}) in the flux for the first order model and in the source term for the second order one. This is why we apply explicit time integration method. Denoting the time step by $\Delta t$, the density at the time step $t^{n+1}$ is obtained by using an explicit Euler method with the splitting technique between the transport and the source terms
\begin{equation}\label{eq:scheme}
      \left\{
           \begin{array}{l}
              \ds{U_{i}^{*}=U_{i}^{n}-\frac{\Delta t}{|C_{i}|}\sum_{j=1}^{N_{i}}|e_{ij}|\mathcal{F}(U^{n}_{i},U^{n}_{j},\vec{n}_{ij})},\\
              \ds{U_{i}^{n+1}=U_{i}^{*}+\Delta t S(U_{i}^{*})},
           \end{array}
      \right.
\end{equation}
where the numerical flux function $\mathcal{F}$ depends explicitly on $\vec{\mu}^{n}$. The stability is achieved under the CFL condition $\displaystyle{\Delta t\leq \alpha \cdot\min_{i=1,...,N}\textrm{diameter}(C_{i})/\sigma_{i}}$. In case of the second order model $\ds{\sigma_{i}=\max_{j=1,...,N_{i}}(\vec{v}_{i}\cdot\vec{n}_{ij})+\sqrt{P'(\rho_{i})}}$ is the maximal value of the characteristic wave speed of the homogeneous part of the system (\ref{eq:mainSystem}). For the first order model, using the speed-density relation (\ref{eq:speed}), the maximal wave speed is given by $\displaystyle{\sigma_{i}=\max_{i=1,...N}\left|\frac{d}{d\rho}F(\rho)\right|=v_{\max}}$ due to the same argument $|\vec{\mu}|=1$ used for the computation of the coefficient $\xi$ in the Lax-Friedrichs flux (\ref{eq:LaxFriedrichsFlux}). The value of $\alpha$ is set to $0.9$ in the following computations.

\subsection{Solution to the eikonal equation and gradient}\label{sec:eikonalSolver}
To obtain the solution at time step $t^{n+1}$ we need to compute the direction vector $\vec{\mu}$ defined by (\ref{eq:directionVector}). It means that we have to solve the eikonal equation (\ref{eq:eikonal}) and compute the gradient of its solution. Equation (\ref{eq:eikonal}) is a special case of the static Hamilton-Jacobi equation, for which many numerical methods have been developed such as level-set methods \cite{Osher_Sethian, Osher_Fedkiw}, fast marching and fast sweeping methods \cite{Sethian, Tsai_FSM}, semi-lagrangian scheme \cite{Falcone_Ferretti}, finite volume or finite element schemes \cite{Hu_Shu, Bryson_Levy}. We implement the Bornemann and Rasch algorithm \cite{Bornemann_Rasch} belonging to the last of the above approaches thus it is easier to implement on unstructured triangular meshes with respect to other methods. It is a linear, finite element discretization based on the solution to a simplified, localized Dirichlet problem solved by the variational principle.

Having found the potential $\phi$ we calculate its gradient using the nodal $P_{1}$ Galerkin gradient method. It is related to cell $C_{i}$ and is computed by averaging the gradients of all triangles having node ${i}$ as a vertex. In two dimensions it has the form
\begin{equation}\label{eq:gradient}
      \nabla\phi_{i}=\frac{1}{|C_{i}|}\sum_{T_{ij}\in C_{i}}\frac{|T_{ij}|}{3}\sum_{k\in T_{ij}}\phi_{k}\nabla P_{k}|_{T_{ij}},
\end{equation}
where $T_{ij}$ are triangles with the considered node $i$ as a vertex, $k$ counts for vertices of $T_{ij}$ and  
$P_{k}|_{T_{ij}}$ is a $P_{1}$ basis function associated with vertex $k$. 

\subsection{Boundary and initial conditions}
\begin{figure}
\begin{center}
\begin{tabular}{c}
\includegraphics[scale=0.4]{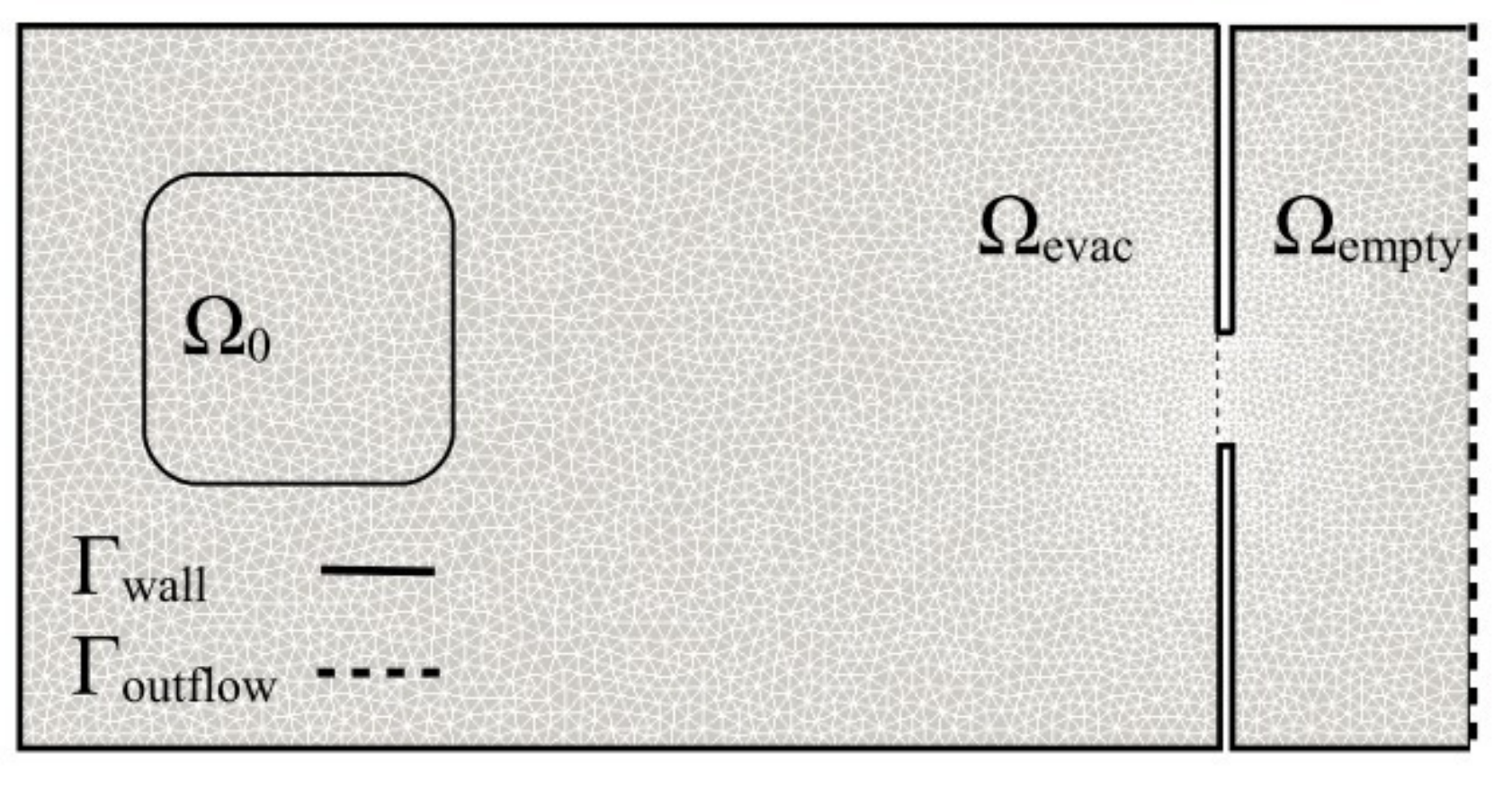}\\
\end{tabular}
\caption{{\emph{\small{Two dimensional domain $\Omega_{evac}\cup\Omega_{empty}=\Omega\subset\mathbb{R}^{2}$  with the boundary $\partial\Omega=\Gamma_{wall}\cap\Gamma_{outflow}$. Initially the density is positive only in the region $\Omega_{0}$}}}}
\label{fig:domain}
\end{center}
\end{figure}

We perform simulations on a two-dimensional domain $\Omega\subset\mathbb{R}^2$ with boundary $\partial\Omega=\Gamma_{\textrm{o}}\cup\Gamma_{\textrm{w}}$ , see Fig~\ref{fig:domain}. We set the outflow boundary $\Gamma_{o}$ far from the exit of the room through which pedestrians go out so that the outflow rate does not influence the flow through the door. We assume pedestrians cannot pass through walls, but can move along them: we impose free-slip boundary conditions 
\begin{equation}\label{eq:slipCondition}
      \vec{v}\cdot\vec{n}=0,\quad\frac{\partial\rho}{\partial n}=0\qquad\textrm{at }\Gamma_{w}.
\end{equation} 

In order to implement \eqref{eq:slipCondition} in the case of the second order model we compute the fluxes $\mathcal{F}(U_{i},U_{g},\vec{n})$ through boundary facets using an interior state $U_{i}$ and a corresponding ghost state $U_{g}$. In particular, we choose
\begin{equation}\label{eq:ghostCell}
      \rho_{g}=\rho_{i},\qquad\vec{v}_{g}=\vec{v}_{i}-2(\vec{v}_{i}\cdot\vec{n})\vec{n}
\end{equation}
at wall boundary $\Gamma_{w}$ and $\displaystyle{\rho_{g}=0.1\rho_{\max}}$, $\displaystyle{\vec{v}_{g}=v_{\max}\vec{n}}$ for the outflow $\Gamma_{o}$. 
\begin{rem}\label{rem:wall}
For the numerical flux function $\mathcal{F}(U_{i},U_{g},\vec{n})$ we use the HLL approximate Riemann solver \cite{HLL}. However, our numerical simulations show that the condition (\ref{eq:slipCondition}) is not satisfied at the wall boundary. Sub iterations would be needed at each time step to converge to the correct solution.To reduce the computational cost, after computing $U^{*}$ in \eqref{eq:scheme} we set to zero at wall boundary nodes the  component of the velocity normal to the boundary. Adding the source term preserves the slip-wall boundary condition.   
\end{rem} 
\begin{rem} 
(Mass conservation) It is essential that there is no loss of the mass through the wall boundary during numerical simulation. The HLL solver \cite{HLL} with the ghost state defined by (\ref{eq:ghostCell}) satisfies this condition when (\ref{eq:slipCondition}) holds. In fact, let us consider four possible combinations of minimum and maximum wave speeds (\ref{eq:waveSpeeds}). Using \eqref{eq:ghostCell} we get $v_{g}^{n}=-v_{i}^{n}=0$ at the wall boundary, where the last equality is due to Remark~\ref{rem:wall}. Then we always have $(\sigma_{i},\sigma_{g})=(-s_{i},s_{i})$, where $s_{i}=\sqrt{P'(\rho_{i})}$. Therefore the flux $\mathcal{F}(U_{i},U_{g},\vec{n})$ is always in the center region of the HLL solver, that is $\sigma_{L}\leq 0\leq\sigma_{R}$ and its first component is zero if $U_{g}$ is defined by (\ref{eq:ghostCell}).
\end{rem}

The first order model consists only of the mass conservation equation and the boundary conditions are imposed by defining directly the fluxes at boundary facets. More precisely, we set $\mathcal{F}=0$ at $\Gamma_{w}$ and $\mathcal{F}=\rho_{\max}V(\rho_{\max})$ at $\Gamma_{o}$.\newline

In our simulation we consider initial conditions of the form $\rho_{0}(x)=\bar{\rho}$ in ${\Omega_{0}}$ and $\vec{v}_{0}(x)=0$, where $\bar{\rho}$ is a positive constant and $\Omega_{0}$ is an area inside the evacuation domain far from the exit, see Fig~\ref{fig:domain}. This means that all pedestrians are placed inside the evacuation room and start to move at $t=0$. 
  
\section{Numerical results}\label{sec:results}

In this section we explore numerically the evacuation dynamics of pedestrians from a room through a narrow exit using  the second order model (\ref{eq:mainSystem}) and the numerical scheme presented in the previous section. We analyze the dependence of its solutions on some of the parameters of the system such as the pressure coefficient $p_{0}$, the adiabatic exponent $\gamma$ and the desired speeds $v_{\max}$. Finally, we study the effect of obstacles on the evacuation from a room through a narrow exit. This analysis is preceded by a numerical error analysis of the scheme and a comparison between the first (\ref{eq:Hughes}) and the second  (\ref{eq:mainSystem}) order model. More details can be found in \cite{InriaReport}.

\begin{table}
\begin{center}
\begin{tabular}{lccc}
\hline
\\
PARAMETER NAME&SYMBOL&VALUE&UNITS
\\
\\
\hline
\\
desired speed&$v_{\max}$&$1-7$&m/s\\
relaxation time&$\tau$&$0.61$&s\\
maximal density&$\rho_{\max}$&$7$&$\textrm{ped/m}^2$\\
pressure coefficient&$p_{0}$&$0.005-10$&$\textrm{ped}^{1-\gamma}\cdot \textrm{m}^{2+\gamma}/\textrm{s}^2$\\
adiabatic exponent&$\gamma$&$2-5$&-\\
density-speed coefficient &$\alpha$ &$7.5$ &-\\
\\
\hline
\end{tabular}
\label{tab:parameters}
\end{center}
\caption{{\emph{\small{Parameters values used in the simulation}}}}
\end{table}
In numerical simulations we use the parameters listed in Table~\ref{tab:parameters}. Maximal velocity $v_{max}$, maximal density $\rho_{max}$ and the response time $\tau$ are chosen from the available literature on experimental studies of pedestrian behaviour, (see \cite{Buchmueller, Seyfried2006}). The values of some of the parameters of model (\ref{eq:mainSystem}), such as $p_{0},\gamma,\alpha$ to authors' knowledge do not have a direct correspondence with the microscopic characteristics  of pedestrian motion. 

We introduce two functionals to analyze the results of simulations: the total mass $\displaystyle{M(t)=\int_{\Omega}\rho(x,t)dx}$ of pedestrians inside the evacuation domain and total evacuation time $\displaystyle{T_{\textrm{evac}}=\int_{0}^{\infty}M(t)dt}$. We use the discrete definitions so the total mass at time step $t^n$ is approximated by $\displaystyle{M^n=\sum_{i=1}^{N}\rho_{i}^{n}|C_{i}|}$ and the discrete total evacuation time by $\displaystyle{T_{\textrm{evac}}=\sum_{n=1}^{\infty}M^{n}\Delta t^{n}}$.

\subsection{Error analysis}
In this section we analyze the accuracy of the numerical scheme presented in Section~\ref{sec:numericalScheme}. More precisely, we first estimate convergence order of the method used to solve the eikonal equation and to compute its gradient. Then we perform the analysis for the fully discrete scheme for the second order model. 
Let $\Omega_{k}$ be a mesh with $N_{k}$ finite volume cells and $\delta_{i}$ be a surface area of the finite volume cell associated with the i-th vertex of the mesh $\Omega_{k}$. We consider the $L^1$ error between the reference solution $u_{ref}$ and the approximated one $u_{h_{k}}$ in the form
\begin{equation}\label{eq:formulaL1error}
      E_{k} = \sum_{i=1}^{N_{k}}|u_{h_{k}}^{i}-u_{ref}^{i}|\delta_{i}.
\end{equation} 
We assume that
\begin{equation}\label{eq:errorI}
      E_{k}^i = u_{h_{k}}^i-u_{ref}^i=Ch_{k}^{p}+h.o.t,\qquad C-constant.
\end{equation}  
The grid-spacing parameter $h_{k}$ has the form $h_{k}=\sqrt{N_{ref}/N_{k}}$, where $N_{ref}$ is the number of finite volume cells for the reference mesh (see \cite{Roache}). When an explicit, analytic solution is available we replace $N_{ref}$ with the area of the domain $|\Omega|$. We estimate the order of convergence $p$ using the least square method applied to the logarithm of the equation  (\ref{eq:errorI})  with neglected higher order terms. 

\subsubsection{Eikonal equation and gradient calculator}
We are going to estimate the order of the algorithm presented in the previous section to solve the eikonal equation and to compute the gradient of its solution.  In the following tests we use the running cost (\ref{eq:runningCost}) with the speed-density relation given by (\ref{eq:speed}) and $v_{\max}=2$ m/s and $\rho_{\max}=7$ $\textrm{ped/m}^2$. Two cases are analyzed.

\begin{figure}[htbp!]
\begin{center}
\begin{tabular}{c}
\includegraphics[scale=0.13]{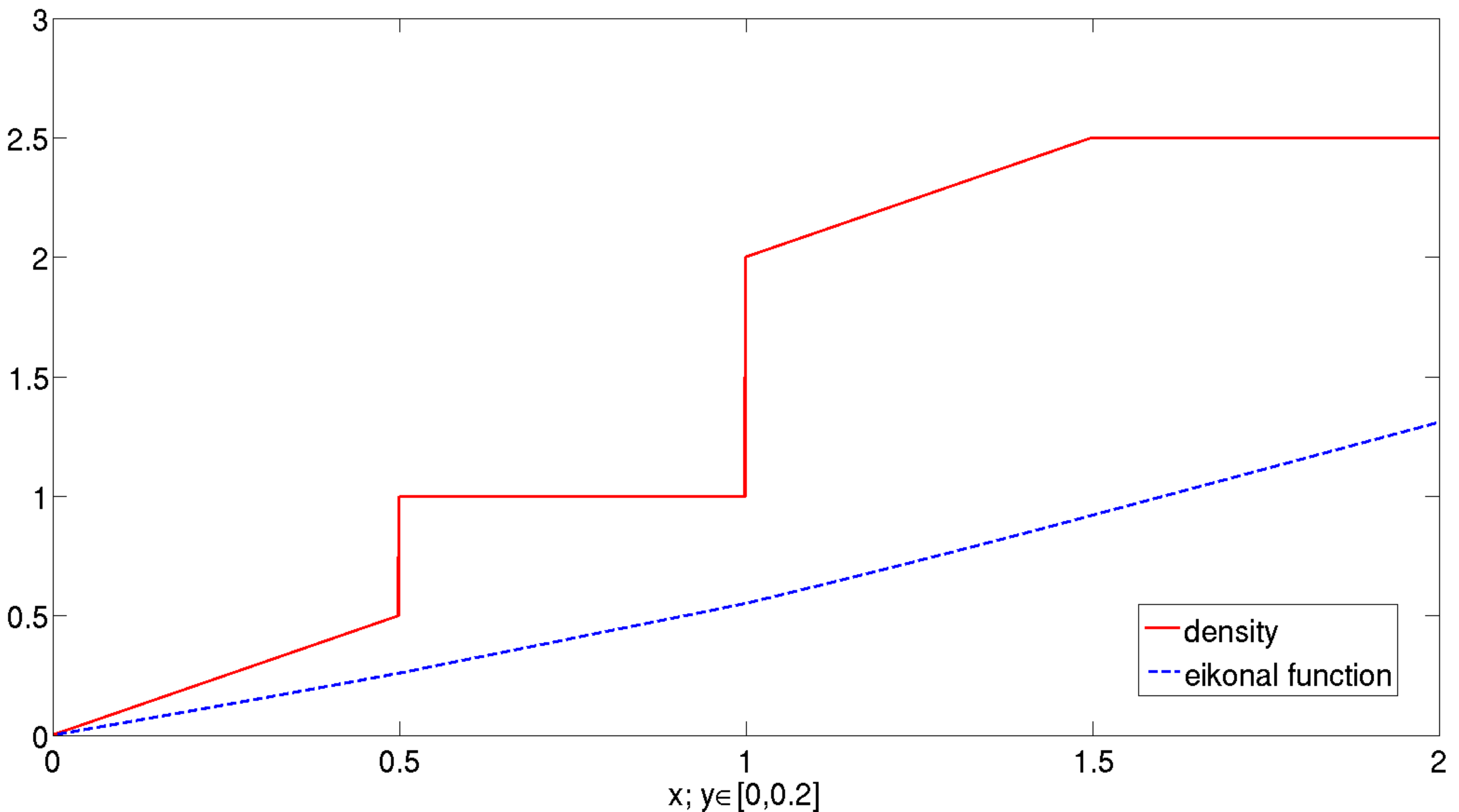}
\end{tabular}
\caption{{\emph{\small{Density distribution (in red) and the explicit solution to the eikonal problem \eqref{eq:eikonal}, \eqref{eq:runningCost} with $\Gamma_{0}=\{(x,y): x=0, y\in[0,0.2]\}$ (in blue) as functions of the x-variable in $\Omega=[0,2\textrm{ m}]\times[0,0.2\textrm{ m}]$ }  }}}
\label{fig:discontDensity}
\end{center}
\end{figure}

\textbf{Test 1:} We consider a domain $\Omega=[0,2\textrm{ m}]\times[0,0.2\textrm{ m}]$ with an outflow localized at $x=0$ and density distribution $\rho(x)=x\cdot\chi_{[0,0.5)}+1\cdot\chi_{[0.5,1)}+(x+1)\cdot\chi_{[1,1.5)}+2.5\cdot\chi_{[1.5,2]}$ $\textrm{ped/m}^2$, where $\chi_{A}$ is zero outside the set $A$, see Fig~\ref{fig:discontDensity}. The solution to the eikonal equation (\ref{eq:eikonal}) and its gradient can be found explicitly from
\begin{displaymath}
      \phi(x)=\int_{0}^{x}\frac{1}{V(\rho)}dx,\qquad\textrm{and}\qquad\vec{\nabla\phi}=\left[\frac{1}{V(\rho(x))},0\right].
\end{displaymath}

\begin{figure}[htbp!]
\begin{center}
\begin{tabular}{c}
\includegraphics[scale=0.26]{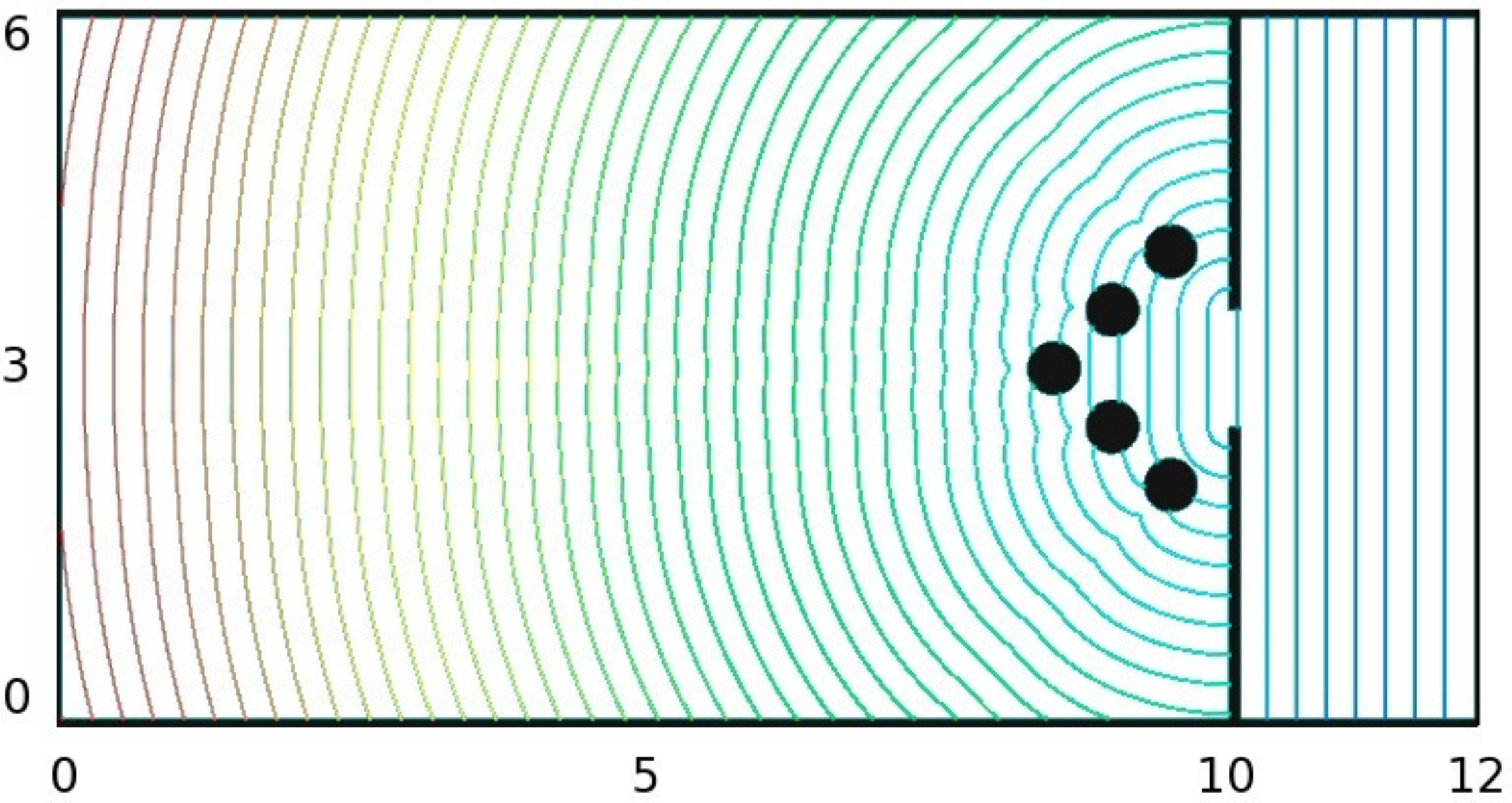}
\end{tabular}
\caption{{\emph{\small{A room $10\textrm{ m}\times 6\textrm{ m}$ with a $1$ m wide, symmetrically placed exit and columns with radius $r=0.23$ m centered at $(9.5,2), (9.0,2.5), (8.5,3), (9.0,3.5), (9.5,4)$.  }}}}
\label{fig:eikonalColumns}
\end{center}
\end{figure}

\textbf{Test 2:} We analyze a domain $\Omega=[0,10\textrm{ m}]\times[0,6\textrm{ m}]$ with an exit of width $L=1$ m centered at $(x,y)=(10,3)$ and five columns of radius $r=0.22$ m in its interior, see Fig~\ref{fig:eikonalColumns}. The density is set to zero everywhere and the reference solution is obtained on a very fine grid with $N=136 507$ finite volume cells. 

\begin{figure}[htbp!]
\begin{center}
\begin{tabular}{cc}
\includegraphics[scale=0.11]{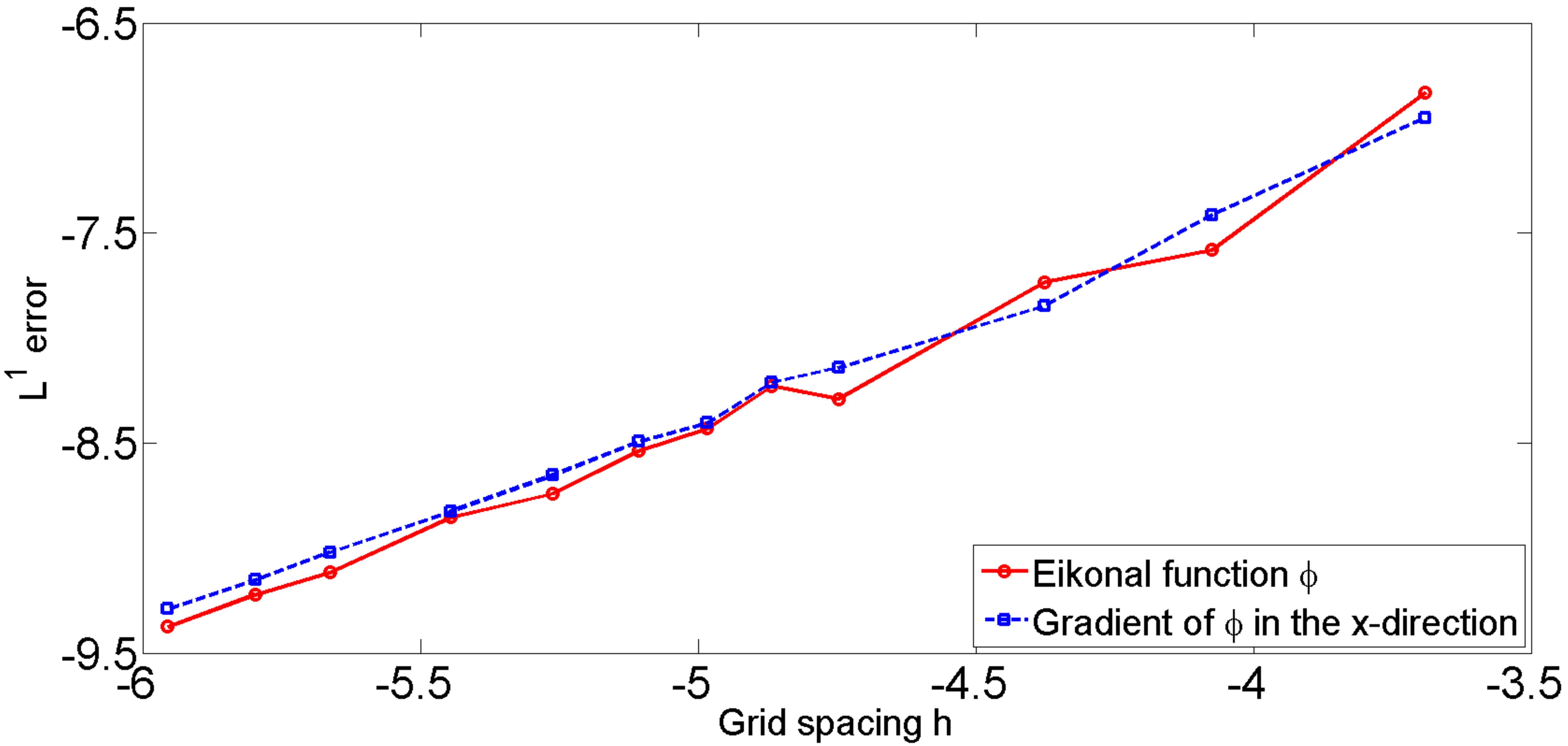}&\includegraphics[scale=0.11]{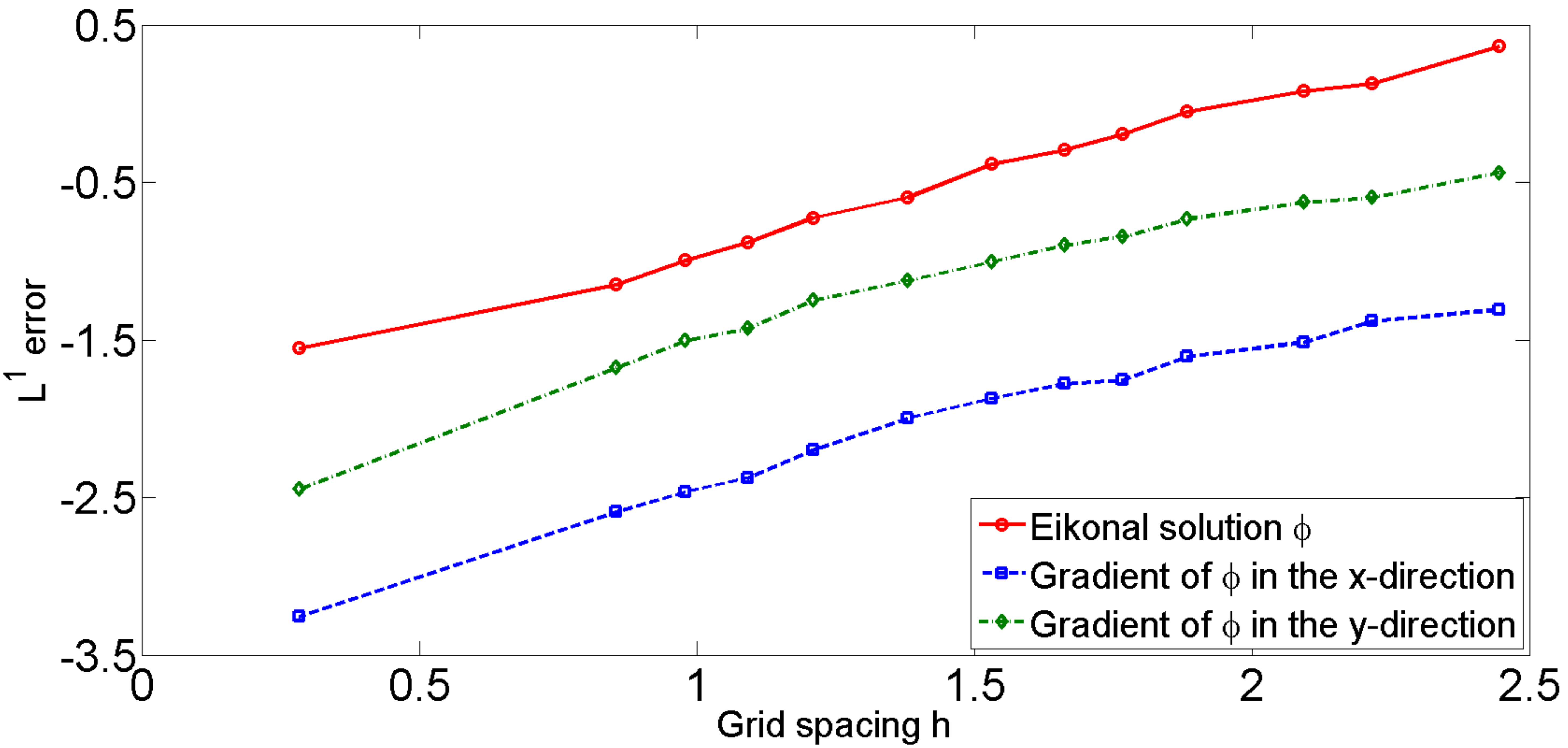}\\
Test 1&Test 2\\
\end{tabular}
\caption{{\emph{\small{$L^{1}$ errors of the eikonal solver and the gradient \eqref{eq:gradient} as functions of the grid-spacing $h$ in the log-log scale for Test 1: a long corridor with a discontinuous density (on the left) and Test 2: a room with five columns (on the right)} }}}
\label{fig:L1errorEikonal}
\end{center}
\end{figure}

In Fig~\ref{fig:L1errorEikonal} we present the loglog plot of the dependence of the $L^1$ error on the grid-spacings $h_{k}$ for the eikonal function and the gradient. In Table~\ref{tab:orderEikonal} we present the corresponding estimates of the convergence order, obtained using the least square method. We observe that the values of the orders are close to one. In particular, gradient is computed with a higher order method and its computation does not decrease the order of the full scheme. 

\begin{table}
\begin{center}
  \begin{displaymath}
  \begin{array}{l|c|c|c}
  \hline
   &\textrm{Eikonal func.}&\textrm{Gradient-x}&\textrm{Gradient-y}\\
  \hline
  \textrm{Test 1}&1.063&1.012&-\\
  \textrm{Test 2}&0.923&0.903&0.881\\
  \hline
  \end{array}
  \end{displaymath}
  \caption{{\emph{\small{ Estimates of the convergence order of the eikonal solver presented in Section~\ref{sec:eikonalSolver} and the method to compute a gradient (\ref{eq:gradient}). }}}}
  \label{tab:orderEikonal}
  \end{center}
  \end{table}

\subsubsection{Second order model}
Now we perform the same error analysis for the fully discrete scheme (\ref{eq:scheme}) for the second order model (\ref{eq:mainSystem}). We consider the same domain as in Test 2 in the previous section, but with initial data $\rho_{0}=1$ $\textrm{ped/m}^2$ in $\Omega_{0}=[1\textrm{ m},5\textrm{ m}]\times[1\textrm{ m},5\textrm{ m}]$, $\vec{v}_{0}=0$ m/s and parameters $v_{\max}=2$ m/s, $\rho_{\max}=7$ $\textrm{ped/m}^2$, $p_{0}=0.005$, $\gamma=2$. 

\begin{figure}[htbp!]
\begin{center}
\begin{tabular}{cc}
\includegraphics[scale=0.11]{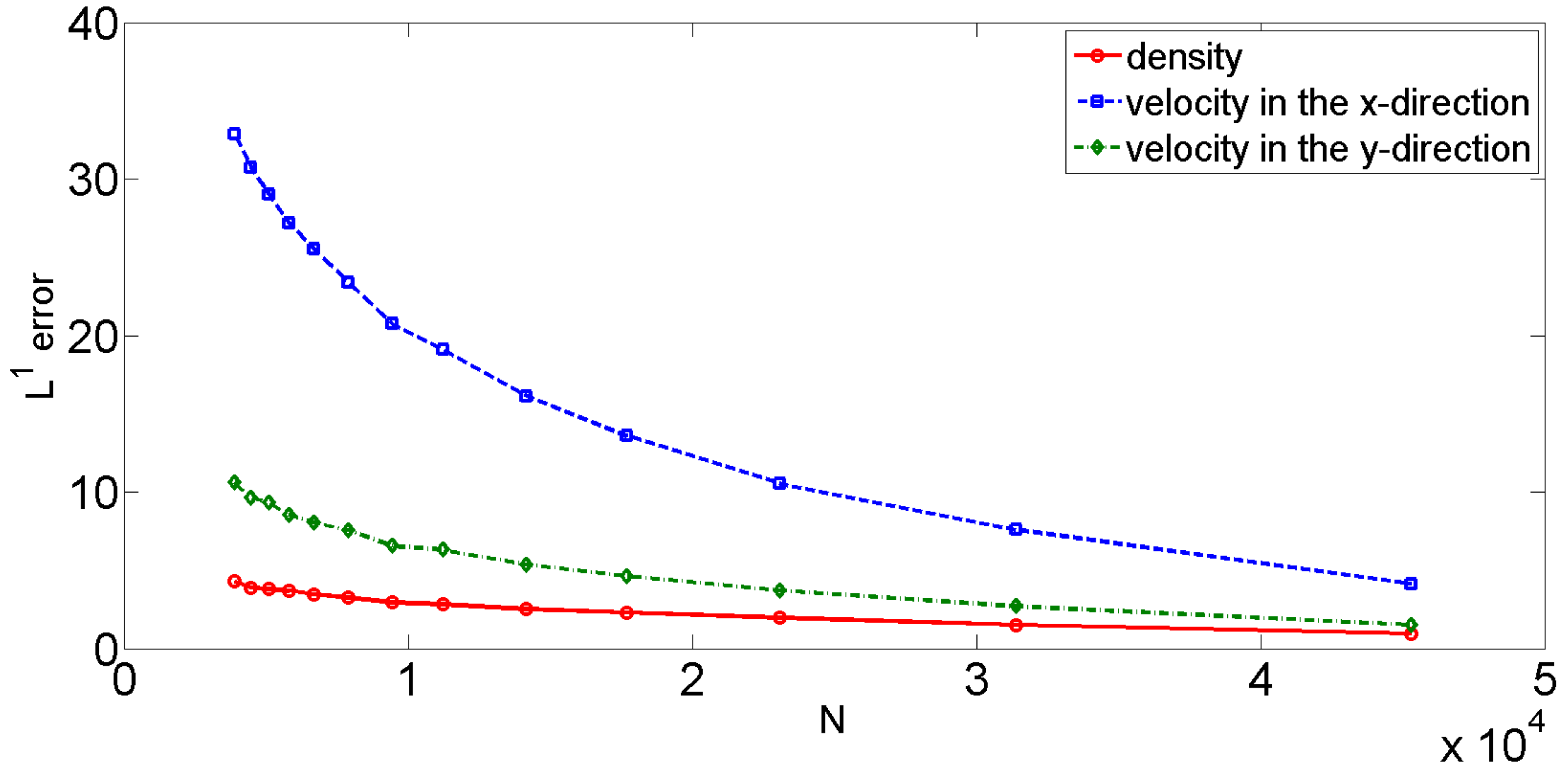}&\includegraphics[scale=0.11]{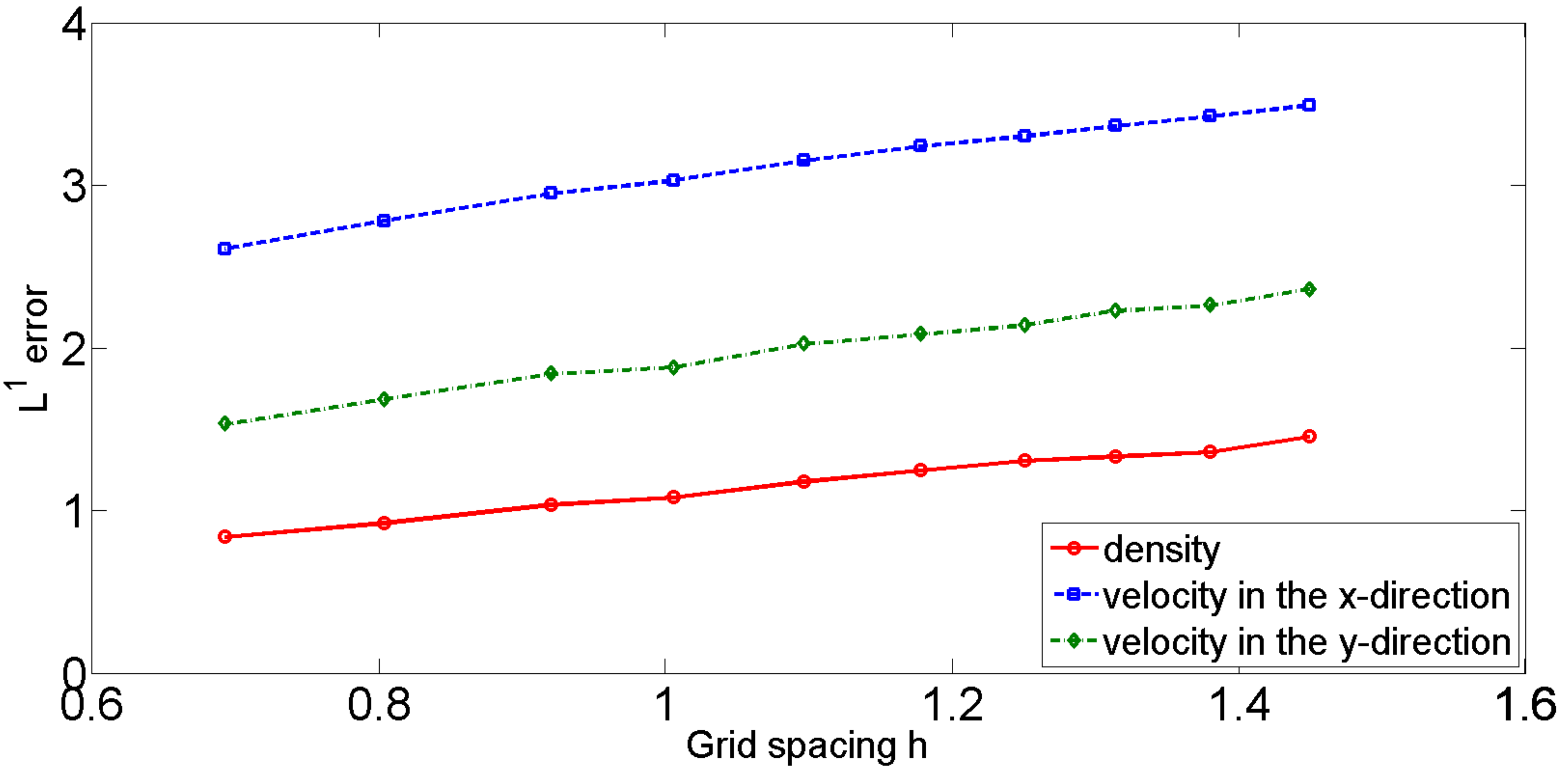}\\
\end{tabular}
\caption{{\emph{\small{$L^1$ errors of the density and the velocity for the second order model \eqref{eq:mainSystem} as a function of the number $N$ of finite volume cells (on the left) and as a function of the grid-spacing h in the log-log scale (on the right).} }}}
\label{fig:L1errorFullSystem}
\end{center}
\end{figure}

The $L^{1}$ error is computed using (\ref{eq:formulaL1error}) with the reference mesh containing $N_{\textrm{ref}}=70 772$ finite volume cells. In Fig~\ref{fig:L1errorFullSystem} we present the dependence at time $t=5$ s of the $L^1$ error on the number of finite volume cells (on the left) and its loglog plot with respect to the grid-spacing parameter $h_{k}$. Least square method gives the following estimates on the order of the full scheme: 0.8 for the density, 1.14 and 1.05 for the velocity in the x and y direction respectively. 

\subsection{Comparison between the first and the second order models}
Now we compare numerically the behaviour of Hughes model \eqref{eq:Hughes} and the second order model \eqref{eq:mainSystem}. We analyze the capability of the models of reproducing the formation of stop-and-go waves and we study the effect on the flow of obstacles placed in the proximity of an exit. 

\subsubsection{Stop and go waves}\label{sec:stopAndGo}

In high density crowd pedestrians experience strong local interactions which can result in macroscopically observed phenomena. One of them, known from vehicular traffic flow, are stop-and-go waves. This corresponds to regions with high density and small speed which propagate backward the flow.  In the case of pedestrians such waves were observed in front of the entrance to the Jamarat Bridge on 12 January 2006 \cite{Helbing_Johansson_Abideen} and studied experimentally for a single lane \cite{Seyfried, Lemercier}. 

\begin{figure}[htbp!]
\begin{center}
\begin{overpic}[scale=0.5,grid,tics=10]{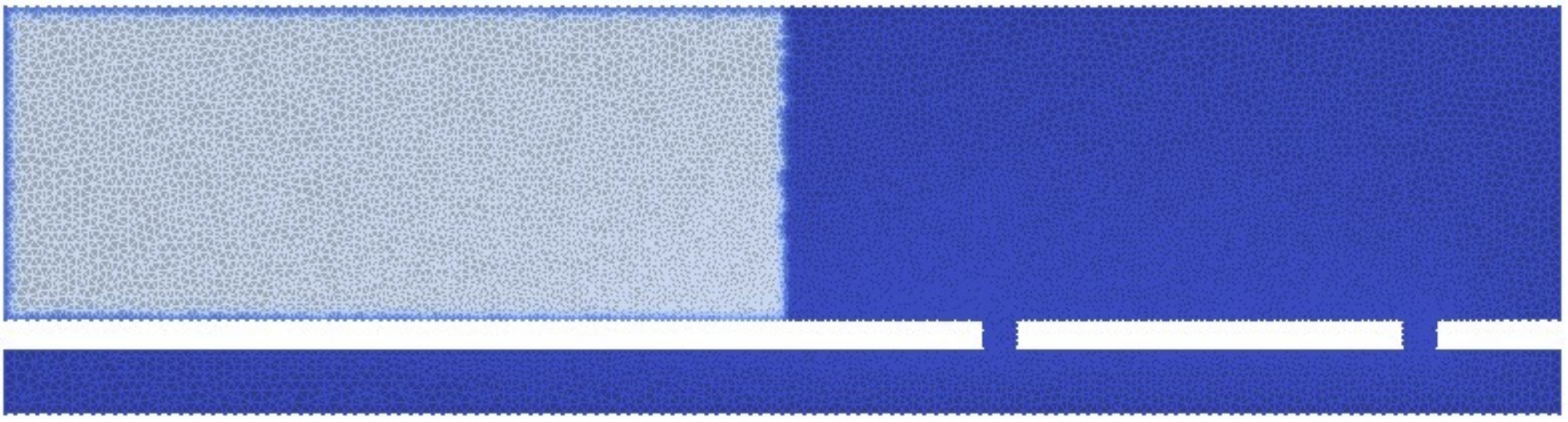}
\end{overpic}
 \caption{{\emph{\small{ A corridor $100\textrm{ m}\times 20\textrm{ m}$ with two, $1.2$ m wide exits centered symmetrically at $x=67$ m and $x=93$ m. }}}}
\label{fig:twoExits_domain}
\end{center}
\end{figure}

\begin{figure}[htbp!]
\begin{center}
\begin{tabular}{lccc}
I&\includegraphics[scale=0.205]{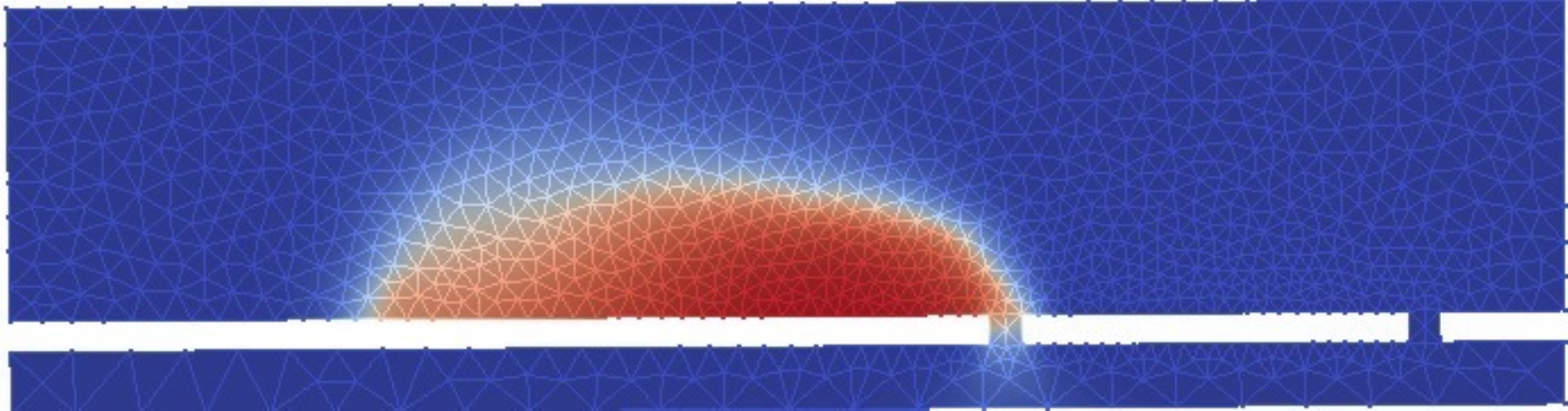}&\includegraphics[scale=0.205]{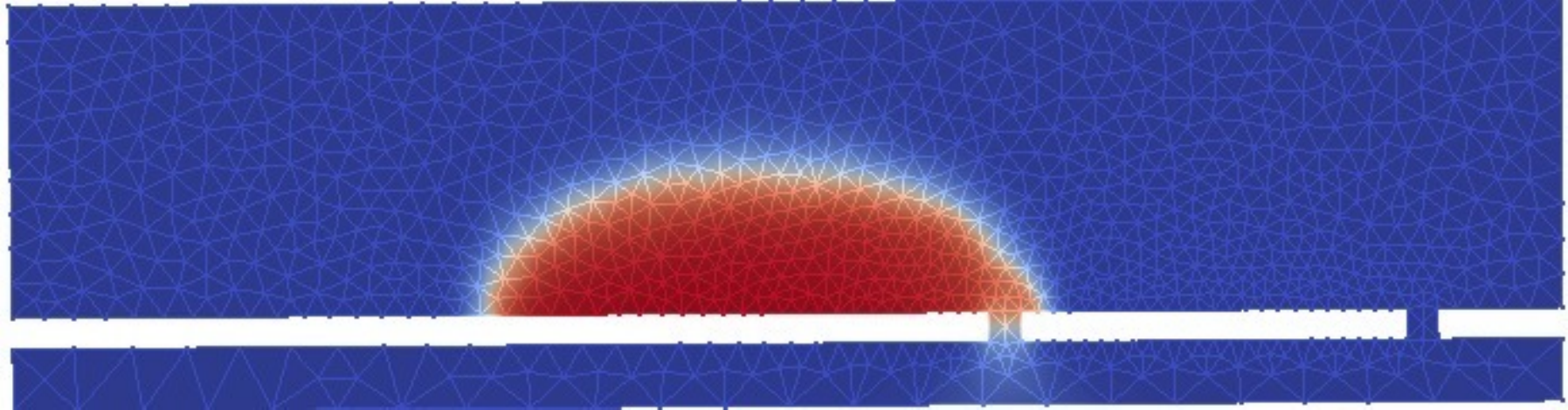}&\includegraphics[scale=0.205]{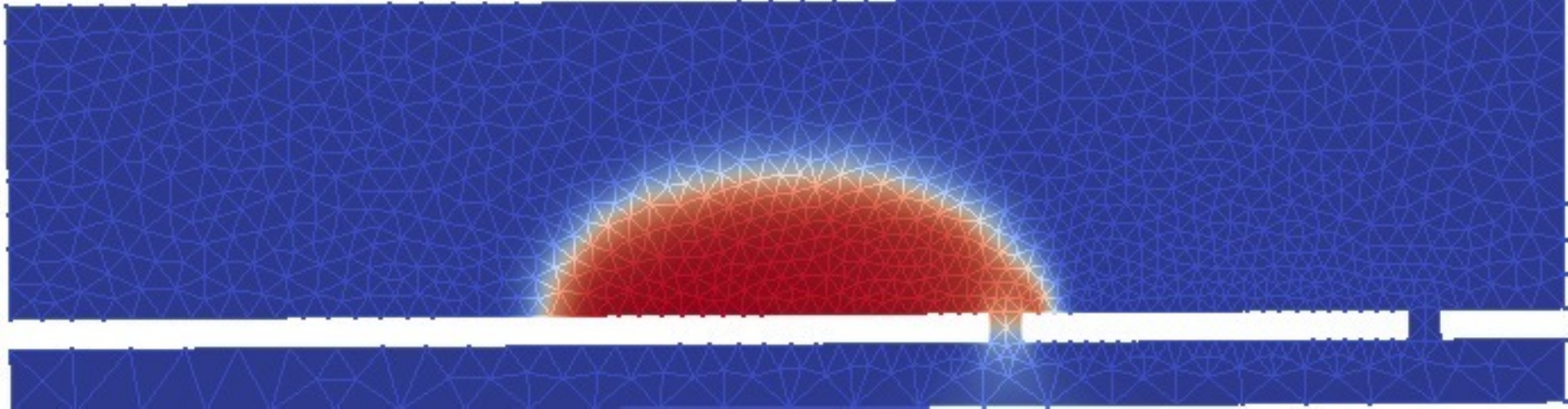}\\
II&\includegraphics[scale=0.18]{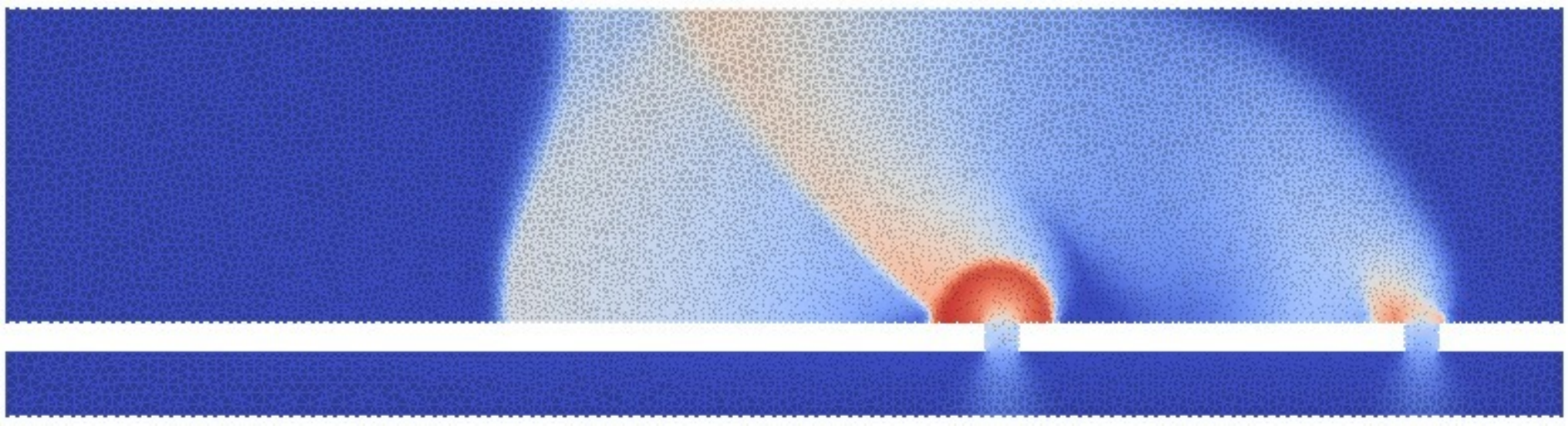}&\includegraphics[scale=0.18]{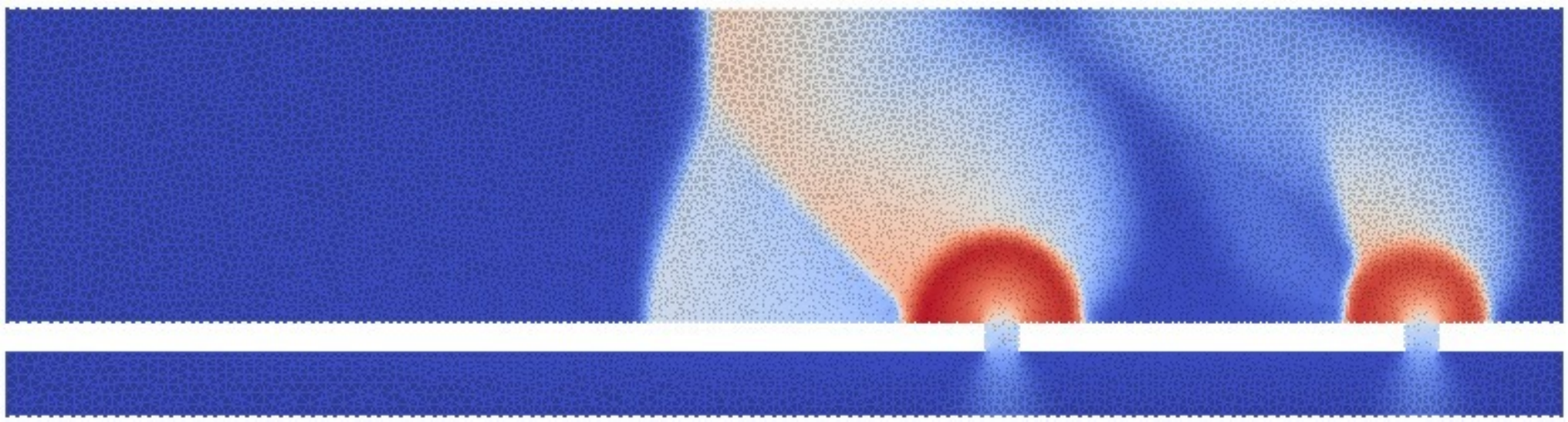}&\includegraphics[scale=0.18]{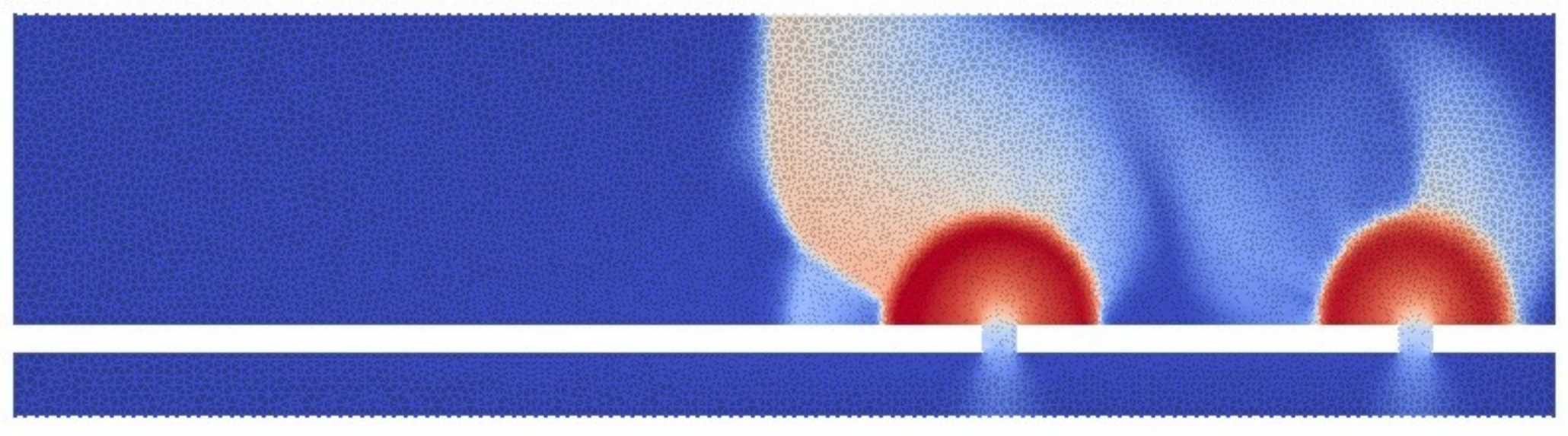}\\
III&\includegraphics[scale=0.18]{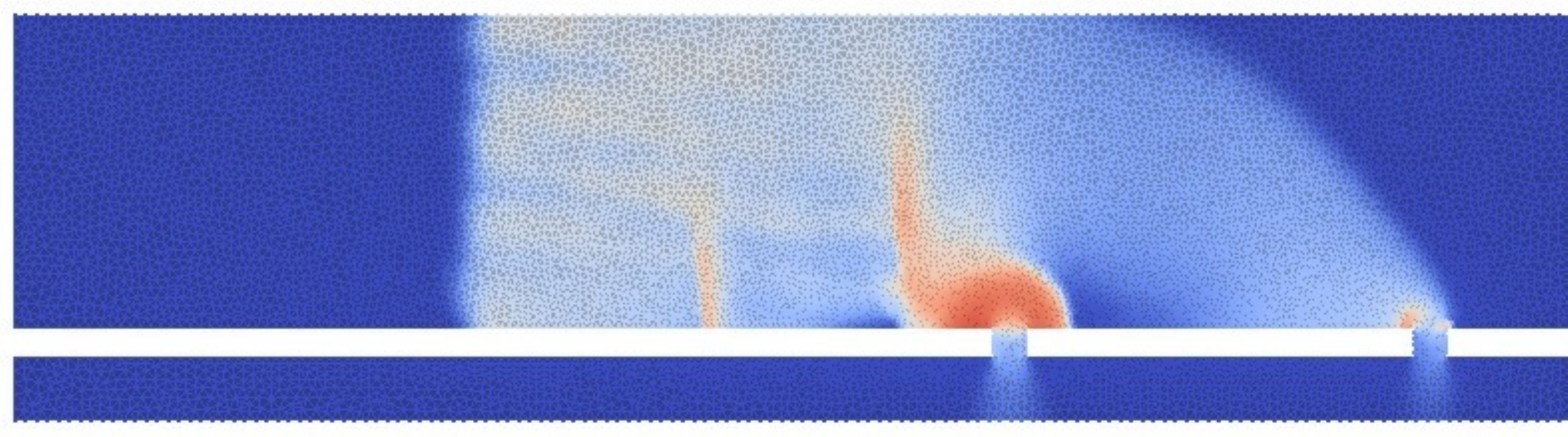}&\includegraphics[scale=0.18]{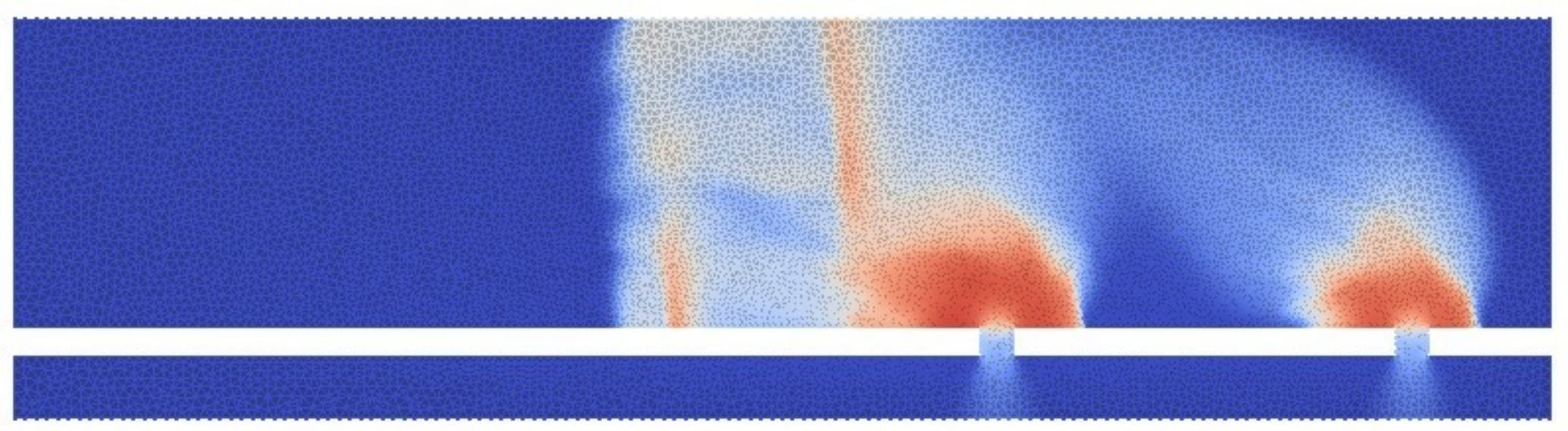}&\includegraphics[scale=0.18]{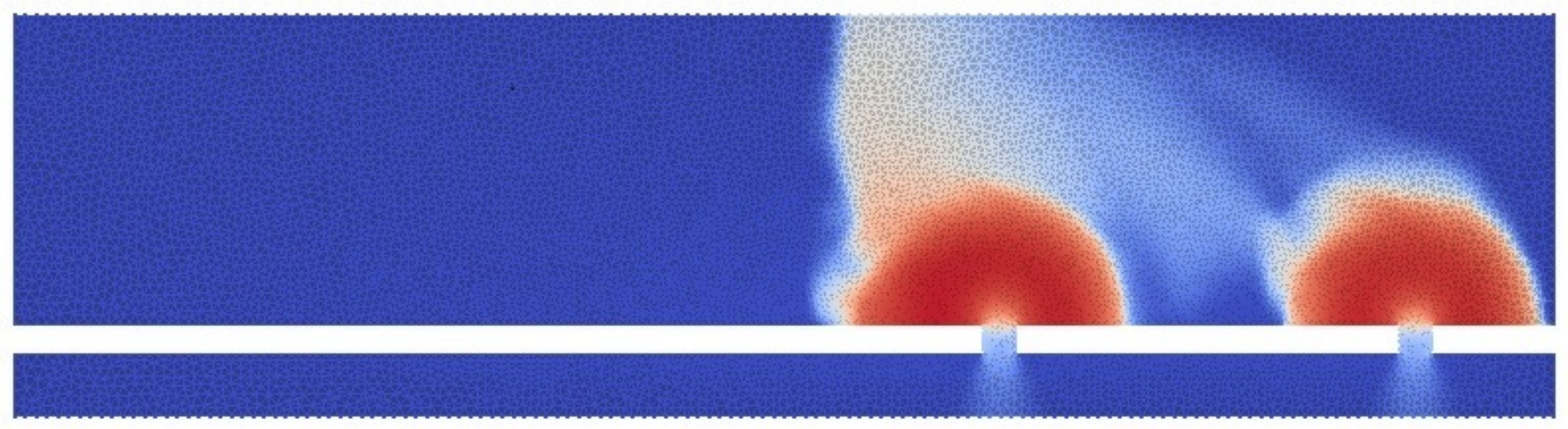}\\
IV&\includegraphics[scale=0.18]{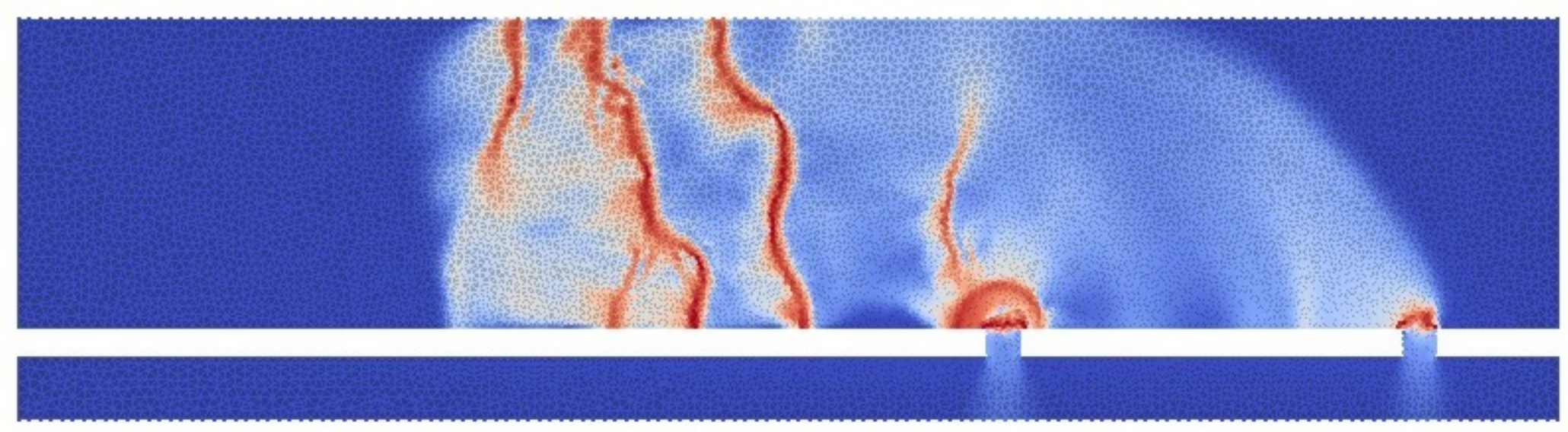}&\includegraphics[scale=0.18]{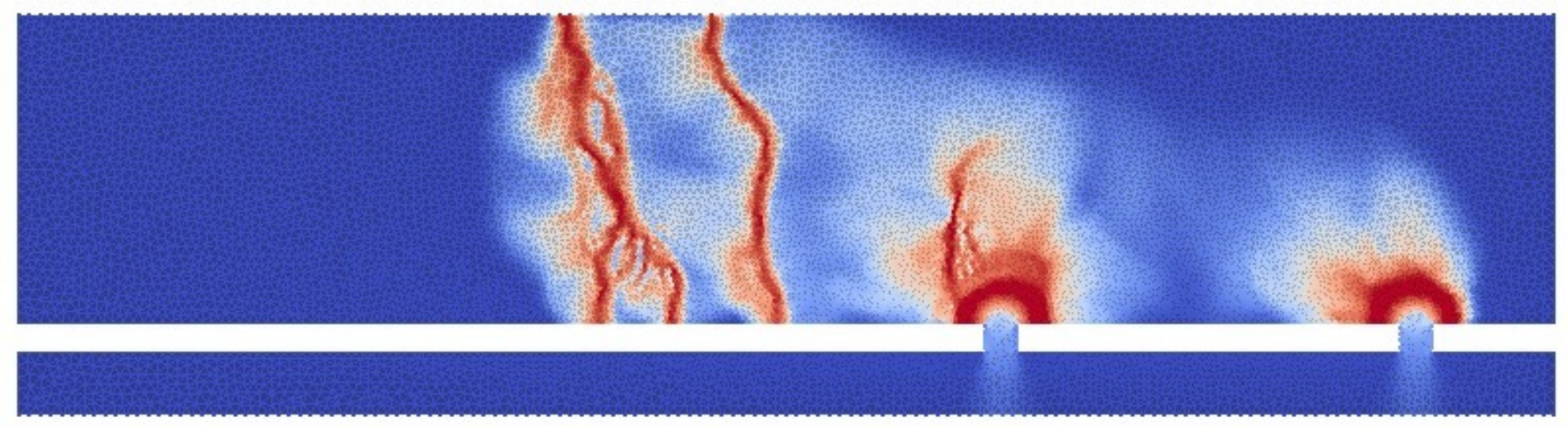}&\includegraphics[scale=0.18]{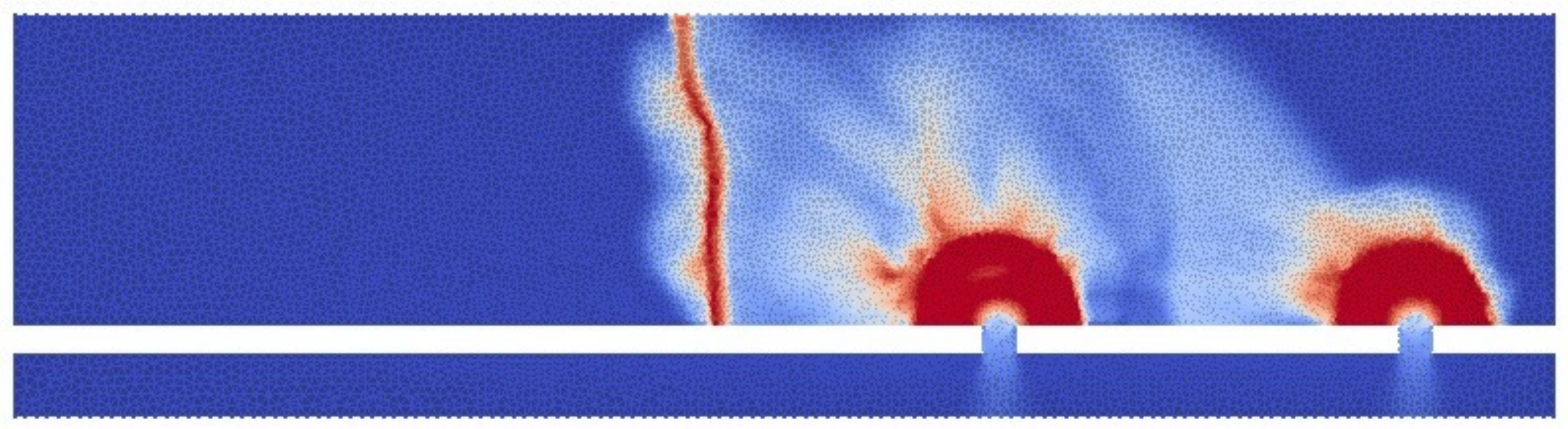}
\end{tabular}
\includegraphics[scale=0.4]{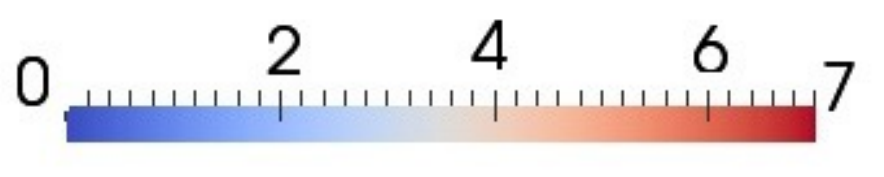}
\caption{{\emph{\small{Density profiles at different times t=30 s (left), t=40 s (middle), t=60 s (right). I: First order model \eqref{eq:Hughes} with $c(\rho)=1\slash v_{\max}$ (Simple model) II: First order model \eqref{eq:Hughes} with $c(\rho)=1\slash V(\rho)$ (Hughes model) III: Second order model with  $c(\rho)=1\slash V(\rho)$ and $p_{0}=0.1$ IV: Second order model with  $c(\rho)=1\slash V(\rho)$ and $p_{0}=0.005$. Other parameters are $v_{\max}=2$ m/s, $\rho_{\max}=7$ $\textrm{ped/m}^2$, $\gamma=2$ and initial data are $\rho_{0}=3$ $\textrm{ped/m}^2$ in $\Omega_{0}=[0,50\textrm{ m}]\times[6\textrm{ m},26\textrm{ m}]$ and $\vec{v}_{0}=0$ m/s.}}}}
\label{fig:twoExits}
\end{center}
\end{figure}

\begin{figure}[htbp!]
\begin{center}
\begin{tabular}{c}
\includegraphics[scale=0.3]{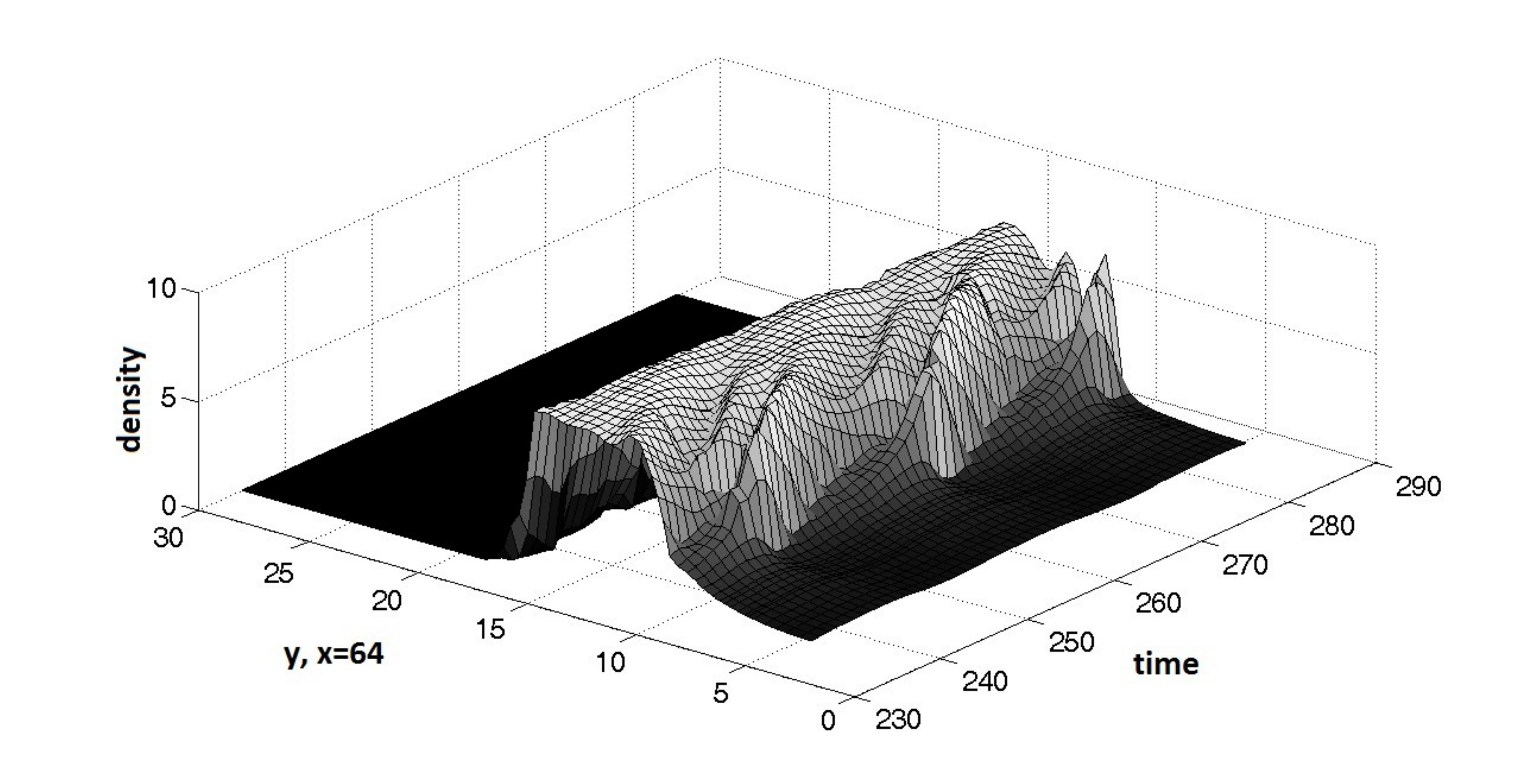}
\end{tabular}
\caption{{\emph{\small{Time evolution of a density profile along a line parallel to the y-axis and passing the center of the first exit $(x,y)=(64,6)$.}}}}
\label{fig:desntiyStopAndGo}
\end{center}
\end{figure}

Simulating stop-and-go waves can be one of the criteria to validate a model, see for example \cite{Lemercier} in case of microscopic description. In order to verify if macroscopic models are able to reproduce this phenomenon we consider a corridor $100\textrm{ m}\times 20\textrm{ m}$ with two, $1.2$ m wide exits centered symmetrically at $x=67$ m and $x=93$ m, see Fig~\ref{fig:twoExits_domain}. 
The initial density is $\rho_{0}=3$ $\textrm{ped/m}^2$ in $\Omega_{0}=[0,50\textrm{ m}]\times[6\textrm{ m},26\textrm{ m}]$ and the initial velocity is set to zero.  Fig~\ref{fig:twoExits} shows the density distribution at different times $t=30, 40, 60$ s in four cases: the first order model \eqref{eq:Hughes} with the running cost function $c(\rho)=1\slash v_{\max}$ (I) and $c(\rho)=1\slash V(\rho)$ (II), the second order model  \eqref{eq:mainSystem} with $c(\rho)=1\slash V(\rho)$ and $p_{0}=0.1$ (III), $p_{0}=0.005$ (IV). Other parameters are $v_{\max}=2$ m/s, $\rho_{\max}=10$ $\textrm{ped/m}^2$, $\gamma=2$. We observe that in case (I) the alternative, further exit is not used by pedestrians, who choose their route only on the basis of the shortest distance, not the shortest time, to the target. The distribution of density given by the first order model (I, II), apart from the proximity of the exit, is very smooth. A similar behaviour is obtained for the second order model with large pressure coefficient $p_{0}=0.1$ (III).  High internal repealing forces between pedestrians prevent from congestion and formation of significant inhomogeneities in the density distribution. Decreasing $p_{0}$ allows for smaller distances between individuals, which causes stronger interactions. As a result for $p_{0}=0.005$ (IV) we observe sub-domains with much higher density. The locations of high density peaks are moving in the opposite direction to the flow. This propagation is clearly observed near the exit where characteristic arching appears. Fig~\ref{fig:desntiyStopAndGo} shows the density distribution profile along the y-direction originating at the center of one of the exits. Stop-and-go waves start at the exit and move backwards.

\subsubsection{Bottlenecks}\label{sec:effectOfObstacle}
\begin{figure}[htbp!]
\begin{tabular}{ccc}
\includegraphics[scale=0.2]{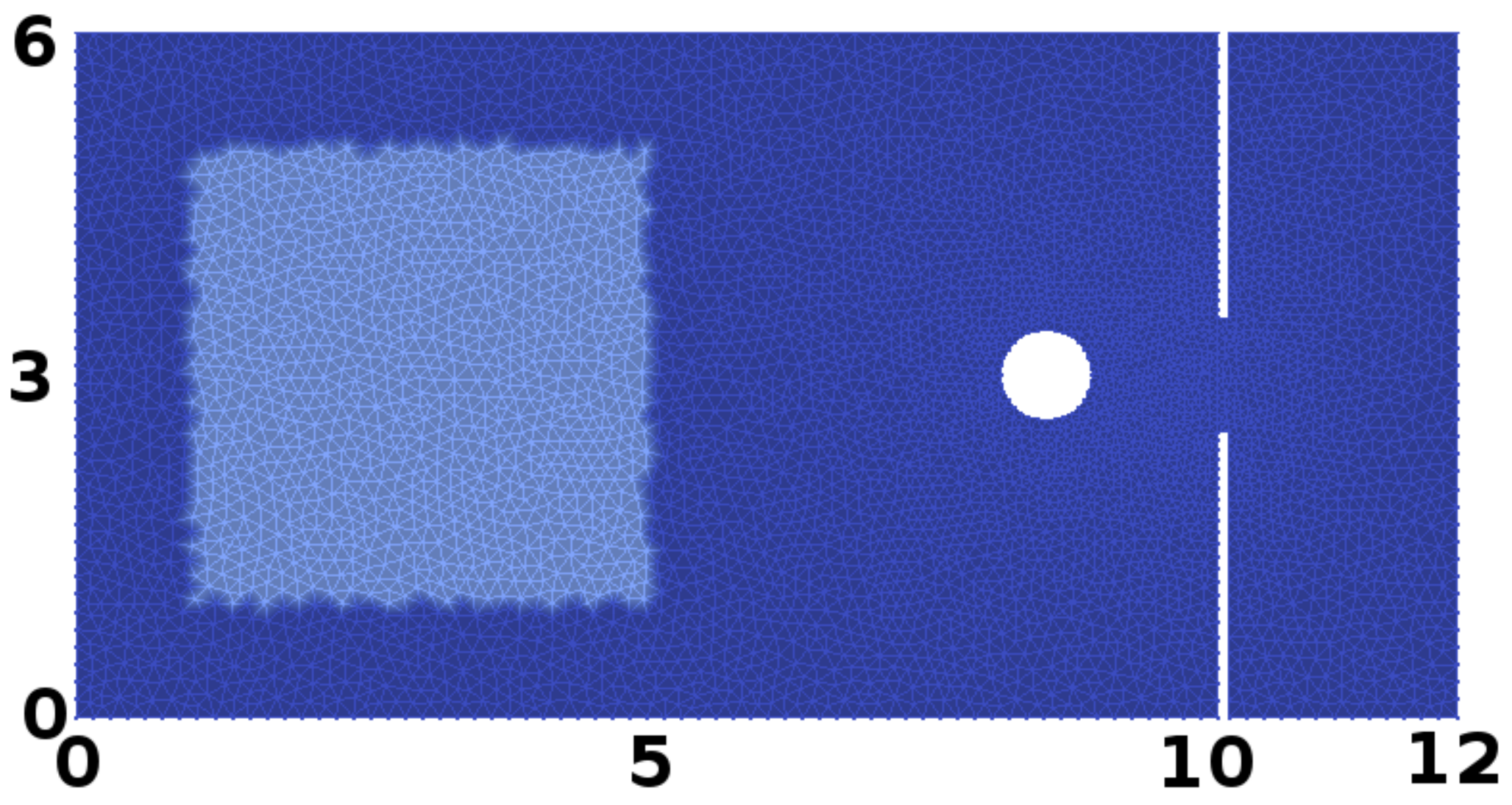}&\includegraphics[scale=0.2]{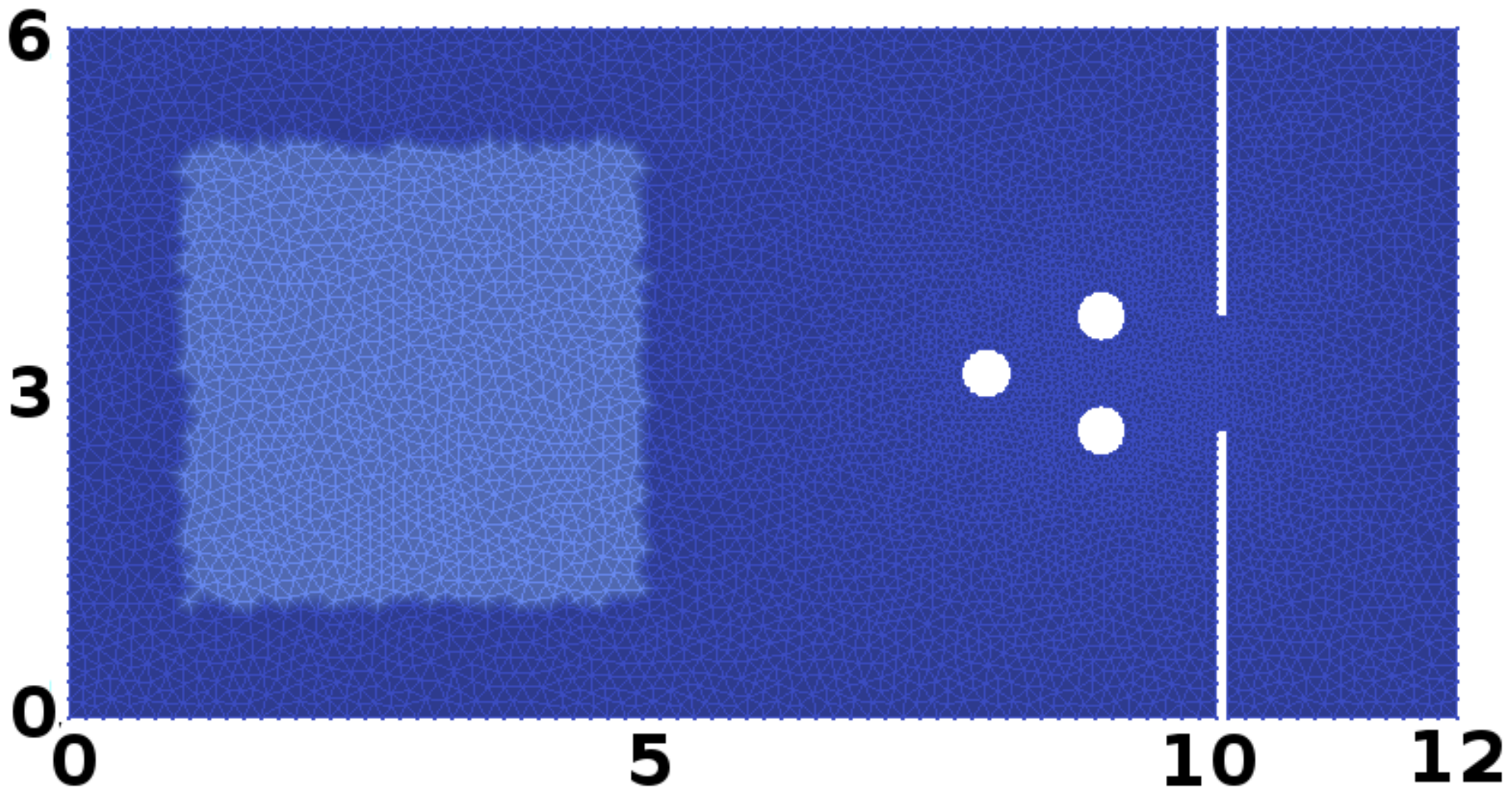}&\includegraphics[scale=0.2]{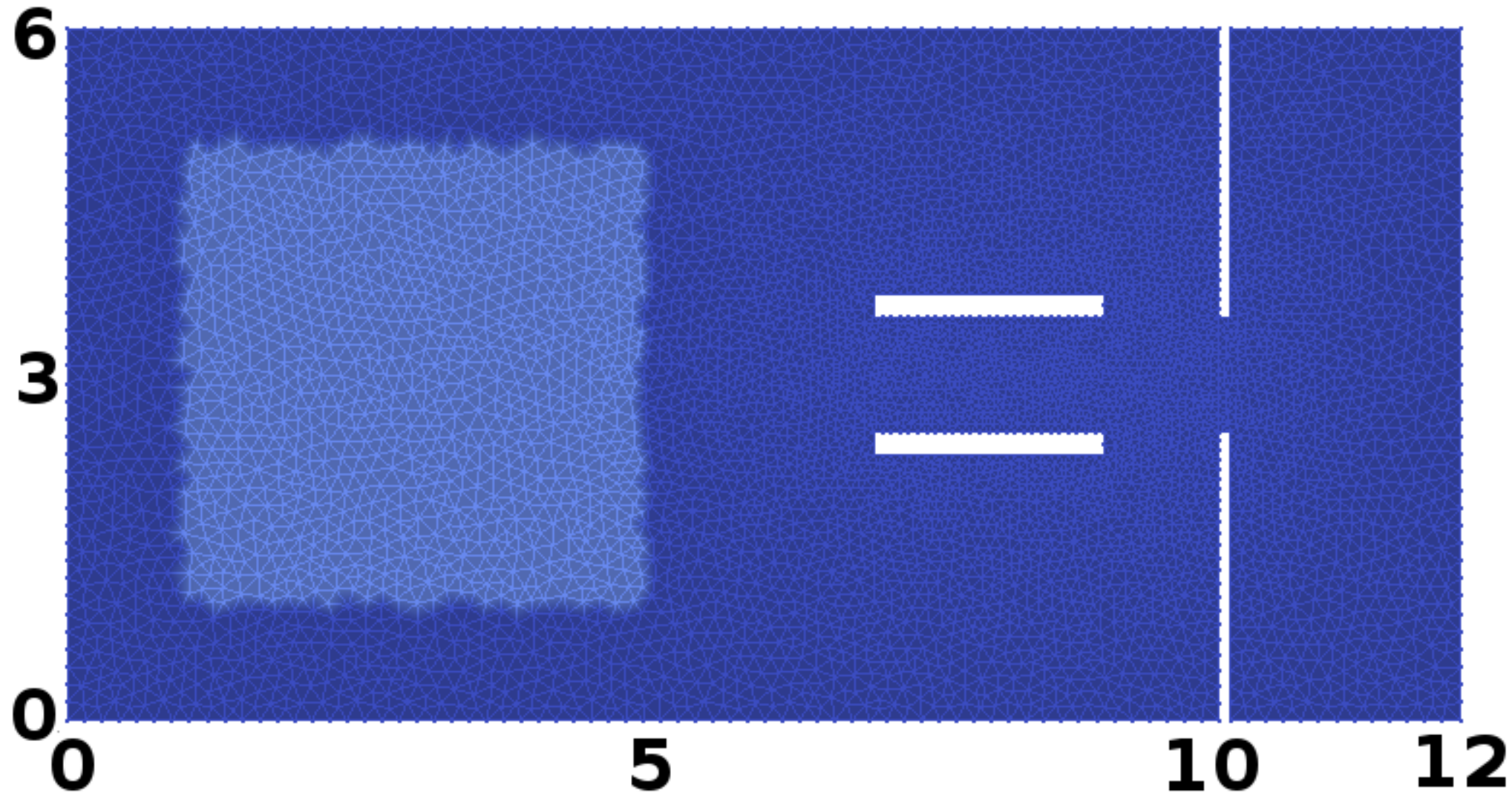}\\
Obstacle 1& Obstacle2& Obstacle3
\end{tabular}
\caption{{\emph{\small{ A room of dimensions $10\textrm{ m}\times 6\textrm{ m}$ with a $1$ m wide, symmetrically placed exit and different obstacles placed in its interior. Obstacle 1: circle centered at $(8.5,3)$ with radius $r=0.3$ m. Obstacle 2: 3 circles centered at $(9,2.5),(8,3),(9,3.5)$ with radius $r=0.2$ m. Obstacle 3:  2 rectangles with $7.5\textrm{ m}\leq x\leq 9\textrm{ m}$ and $2.3\textrm{ m}\leq y\leq 2.5\textrm{ m}$, $3.5\textrm{ m}\leq y\leq 3.7\textrm{ m}$.}}}}
\label{fig:room_obstacles}
\end{figure}

It was experimentally observed that flow of pedestrians through a bottleneck depends on its width (see for example \cite{predtechenskii1978planning, Kretz, HoogendoornDaamen}) and can be significantly slowed down due to clogging at its entrance. Blocking at bottlenecks occurs when the flow of pedestrians towards the door is much higher than the capacity of the exit. The density grows and, as a result, physical interactions between pedestrians increase, slowing down the motion and interrupting the outflow. From the point of view of evacuation strategies it is essential for a mathematical model to be able to capture this effect. We study numerically the behaviour of solutions of the first (\ref{eq:Hughes}) and the second (\ref{eq:mainSystem}) order model during the evacuation through a narrow exit. We consider a room $10\textrm{ m}\times 6\textrm{ m}$ with a $1$ m wide, a symmetrically placed exit and different obstacles placed in its interior, see Fig~\ref{fig:room_obstacles},\newline
\begin{tabular}{ll}
obstacle 1:&one column with radius $r=0.3$ m,\\
obstacle 2:&three columns with radius $r=0.2$ m,\\
obstacle 3:&two walls. 
\end{tabular}

\begin{figure}[htbp!]
\begin{center}
\begin{tabular}{cc}
\includegraphics[scale=0.11]{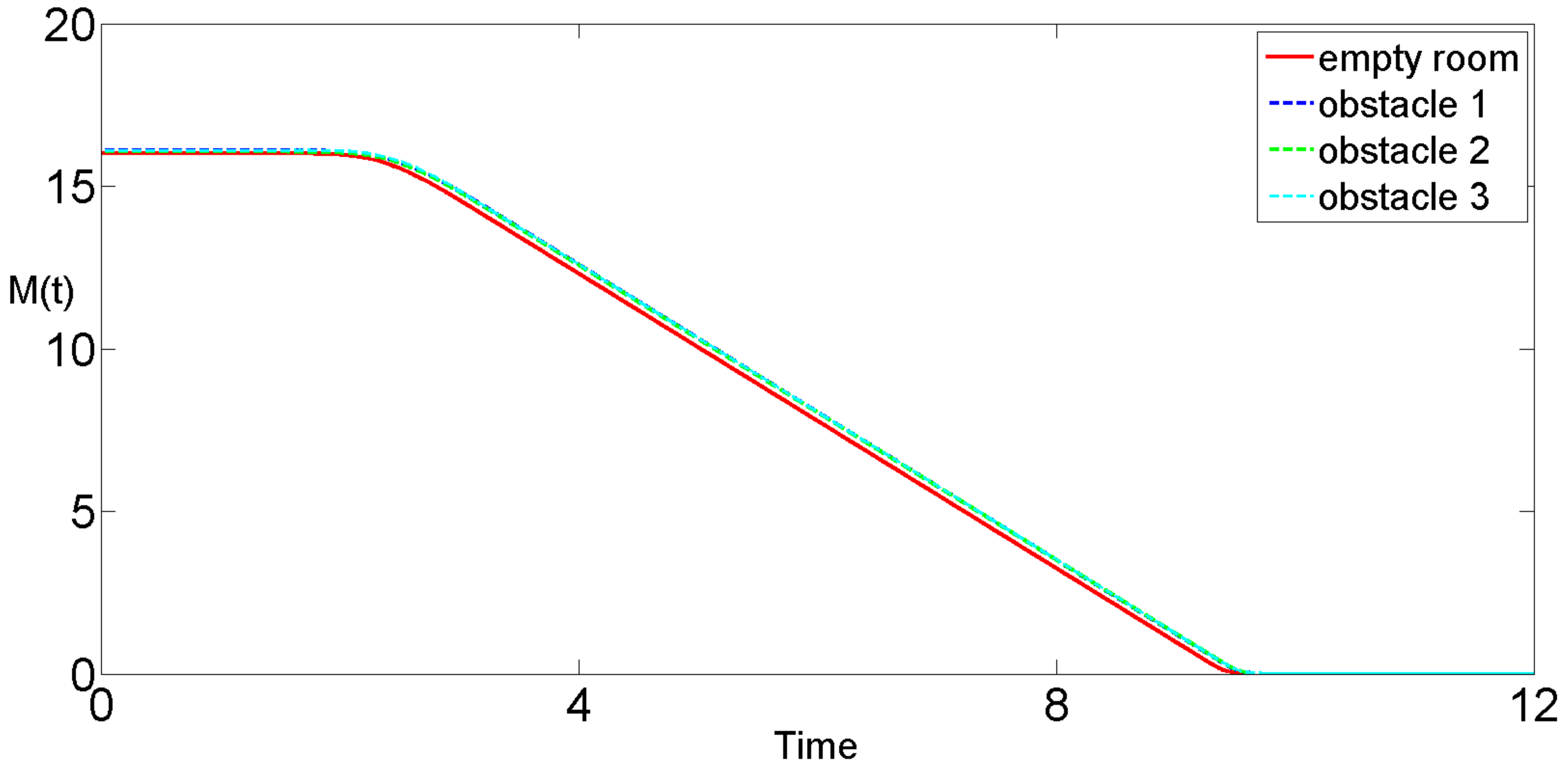}&\includegraphics[scale=0.11]{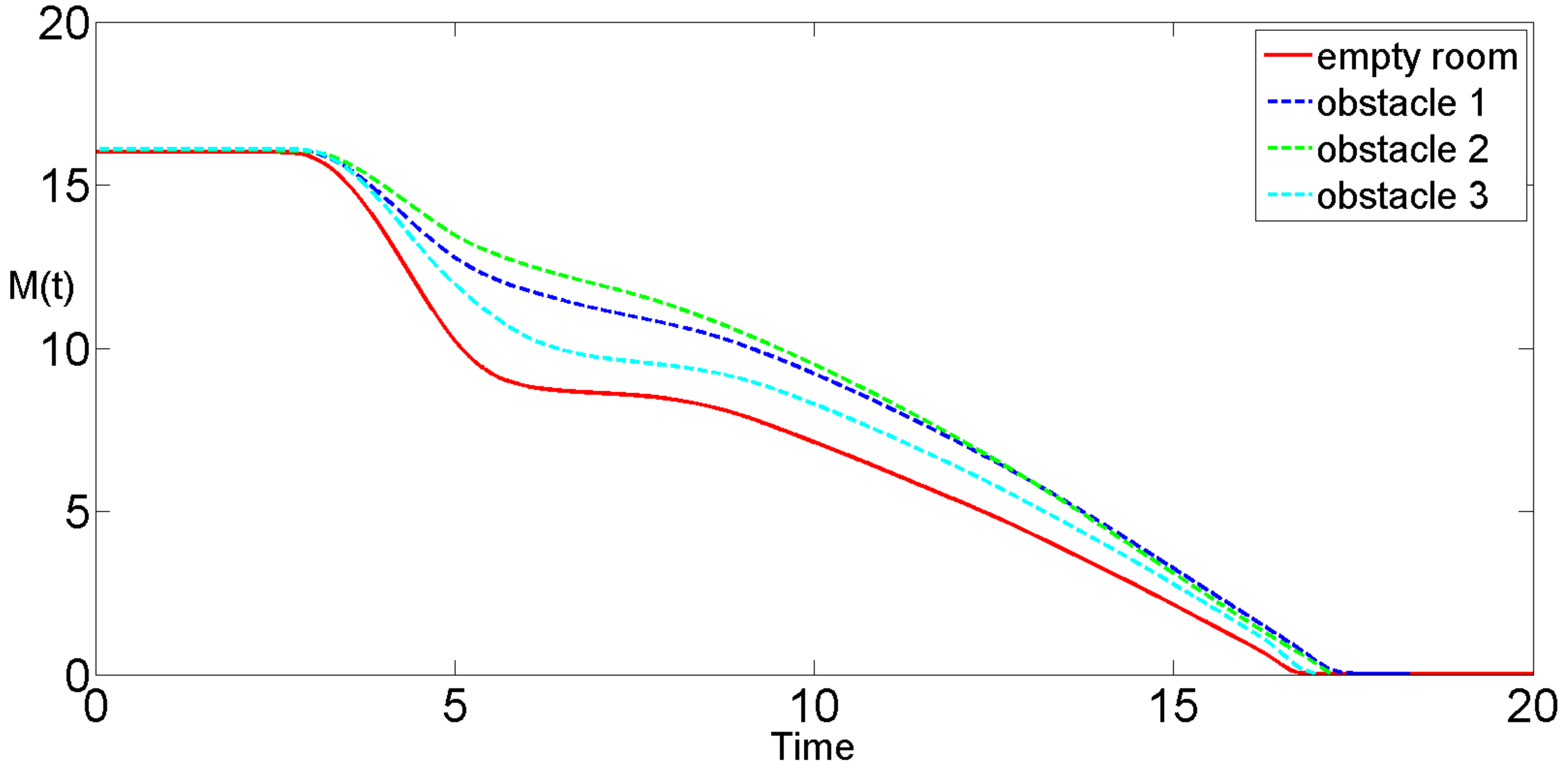}\\
Hughes model&Second order model
\end{tabular}
\caption{{\emph{\small{Time evolution of the total mass $M(t)$ of pedestrians in an empty room and in a room with three different obstacles for the Hughes model \eqref{eq:Hughes} with $v_{\max}=2$ m/s, $\rho_{\max}=7$ $\textrm{ped/m}^2$ (on the left) and for the second order model \eqref{eq:mainSystem} with $v_{\max}=2$ m/s, $\rho_{\max}=7$ $\textrm{ped/m}^2$, $p_{0}=0.005$, $\gamma=2$ (on the right)}. }}}
\label{fig:bottleneck}
\end{center}
\end{figure}

Fig~\ref{fig:bottleneck} presents the time evolution of the total mass of pedestrians remaining inside the room for Hughes' model \eqref{eq:Hughes} (on the left) and for the second order model \eqref{eq:mainSystem} (on the right) with $v_{\max}=2$ m/s, $\rho_{\max}=7$ $\textrm{ped/m}^2$, $P(\rho)=0.005\rho^2$ and with initial data $\rho_{0}=1$ $\textrm{ped/m}^2$ in $\Omega_{0}=[1\textrm{ m},5\textrm{ m}]\times[1\textrm{ m},5\textrm{ m}]$. We observe that in the case of the first order model the total mass $M(t)$ decreases linearly and is the same for all obstacles and the empty room. The outflow is regulated only by the capacity of the door. We do not observe either clogging, which would slow down the decrease of the total mass, or the influence of obstacles. The total evacuation time $T_{evac}$ is basically the same for all situations. On the other hand, in the case of the second order model the empty room experiences a significant decrease of the outflow. It may correspond to clogging at the exit when a sufficiently large number of pedestrians reaches it. Obstacles play a role of a barrier and decrease the flow arriving at the exit. As a result, evacuation is slower but clogging is reduced. 

\subsection{Room evacuation}
We now focus our analysis on the second order model (\ref{eq:mainSystem}) and explore the dependence of its solutions on the system parameters $p_{0},\gamma$ and $v_{\max}$. Our aim is to analyze how the degree of congestion of the crowd and the value of desired walking speed of pedestrians influence the evacuation time.

\subsubsection{Dependence on the parameters $p_{0}$ and $\gamma$}
At first we analyze the dependence on the coefficients of the internal pressure function given by the law for isentropic gases (\ref{eq:pressure}). It describes the repealing forces between pedestrians and prevents from overcrowding. Under emergency and panic conditions the comfort zone of pedestrians, which defines how close they can stay to each other, decreases. As a result, the density of the crowd can increase and impedes the movement leading to discontinuous flow. As we have already seen, in Section (\ref{sec:stopAndGo}), the formation of stop-and-go waves depends on the coefficient $p_{0}$. Now our aim is to study the relation between the strength of the repealing forces between pedestrians and the efficiency of the evacuation.

\begin{figure}[htbp!]
\begin{center}
\begin{tabular}{ll}
\includegraphics[scale=0.11]{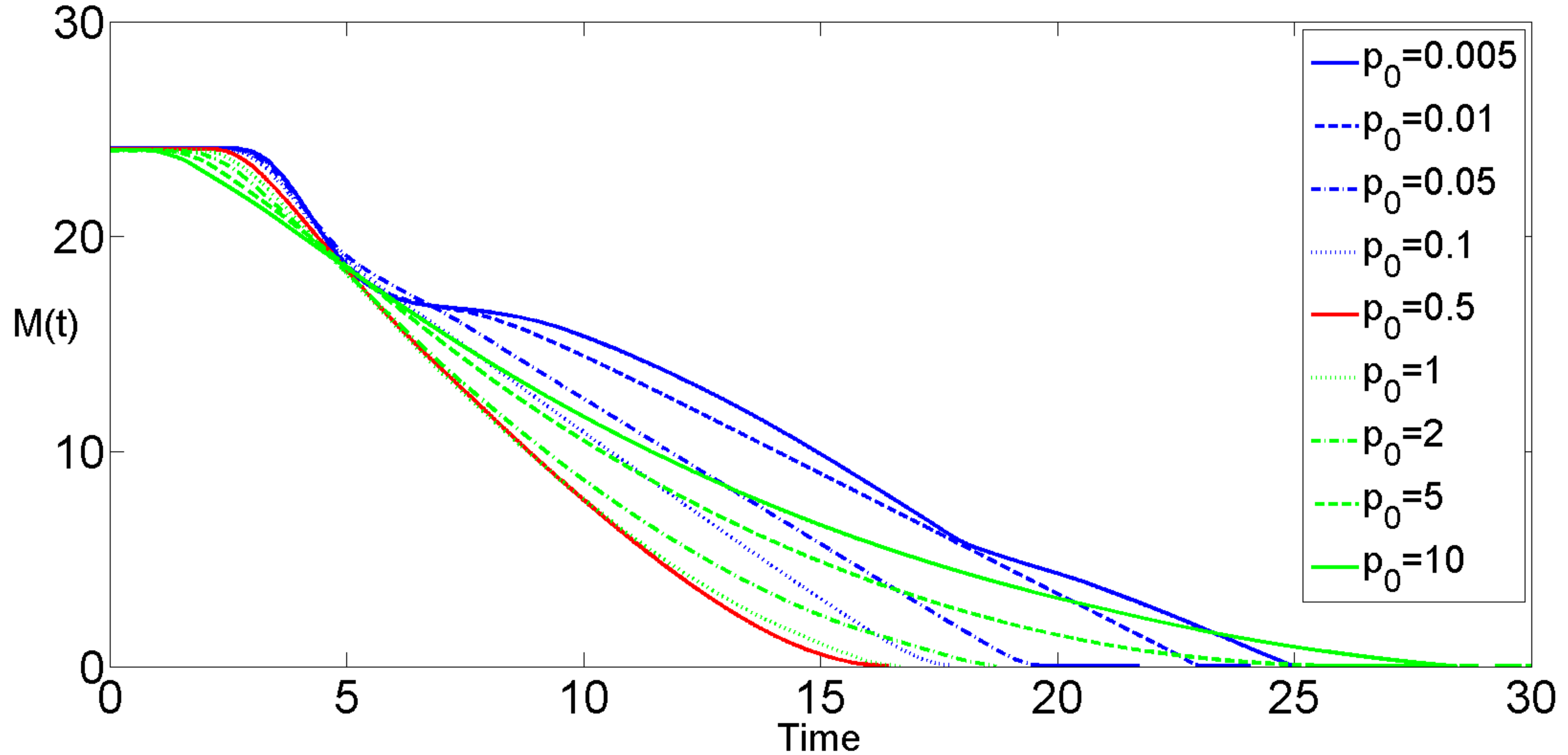}&\includegraphics[scale=0.11]{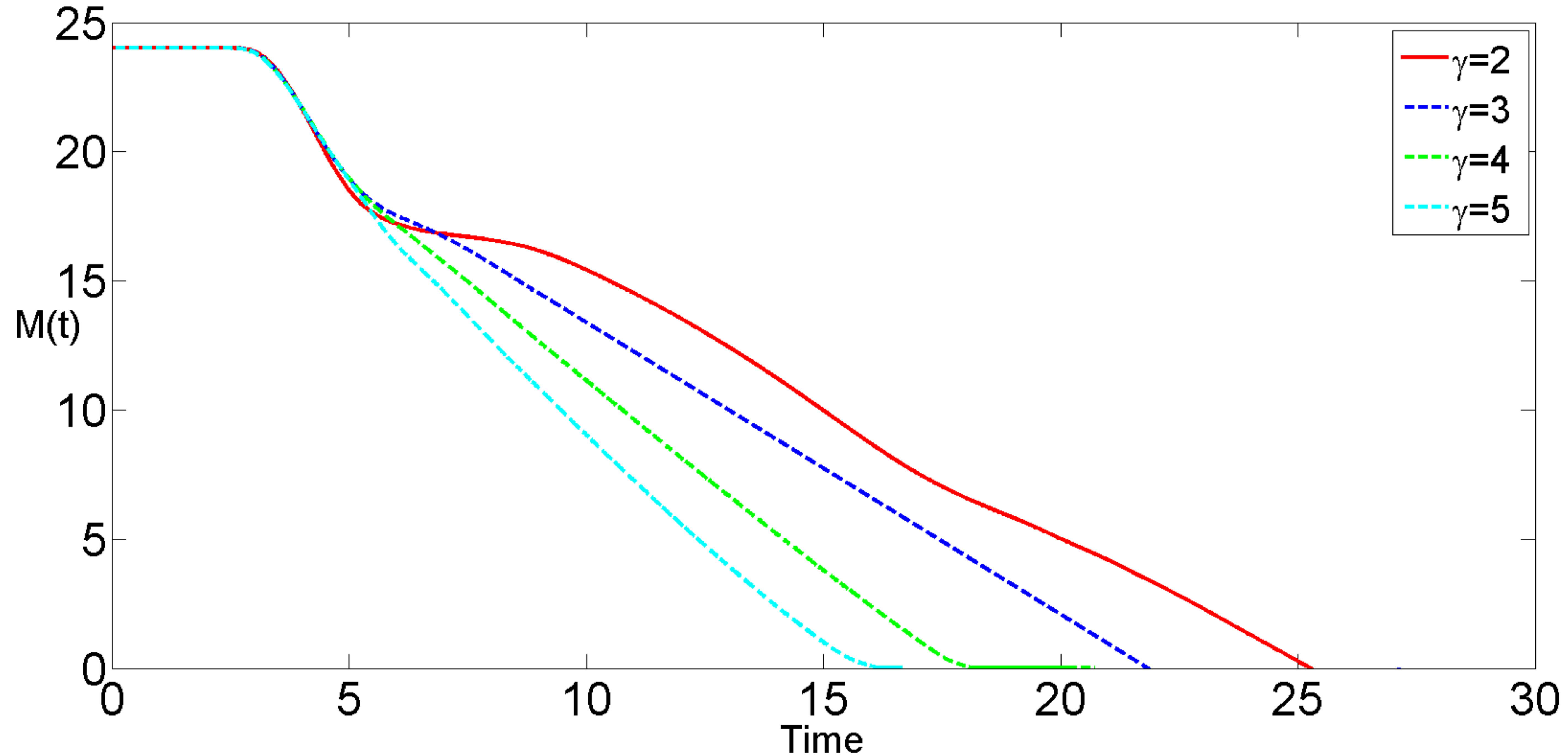}
\end{tabular}
\caption{{\emph{\small{Time evolution of the total mass $M(t)$ of pedestrians inside a room without obstacles for the second order model \eqref{eq:mainSystem} with $v_{\max}=2$ m/s, $\rho_{\max}=7$ $\textrm{ped/m}^2$ $\gamma=2$ for different pressure coefficients $p_{0}=\{5\times 10^{-3},10^{-2},5\times 10^{-2},10^{-1},5\times 10^{-1},1,2,5,10\}$ (on the left) and with $v_{\max}=2$ m/s, $\rho_{\max}=7$ $\textrm{ped/m}^2$, $p_{0}=5\times 10^{-3}$ for different adiabatic exponents $\gamma=2,3,4,5$ (on the right). The initial density is $\rho_{0}=1.5$ $\textrm{ped/m}^2$ in $\Omega_{0}=[1\textrm{ m},5\textrm{ m}]\times[1\textrm{ m},5\textrm{ m}]$, and the initial velocity $\vec{v}_{0}=0$ m/s.}}}}
\label{fig:pressure}
\end{center}
\end{figure}

\begin{figure}[htbp!]
\begin{center}
\begin{tabular}{cc}
\includegraphics[scale=0.15]{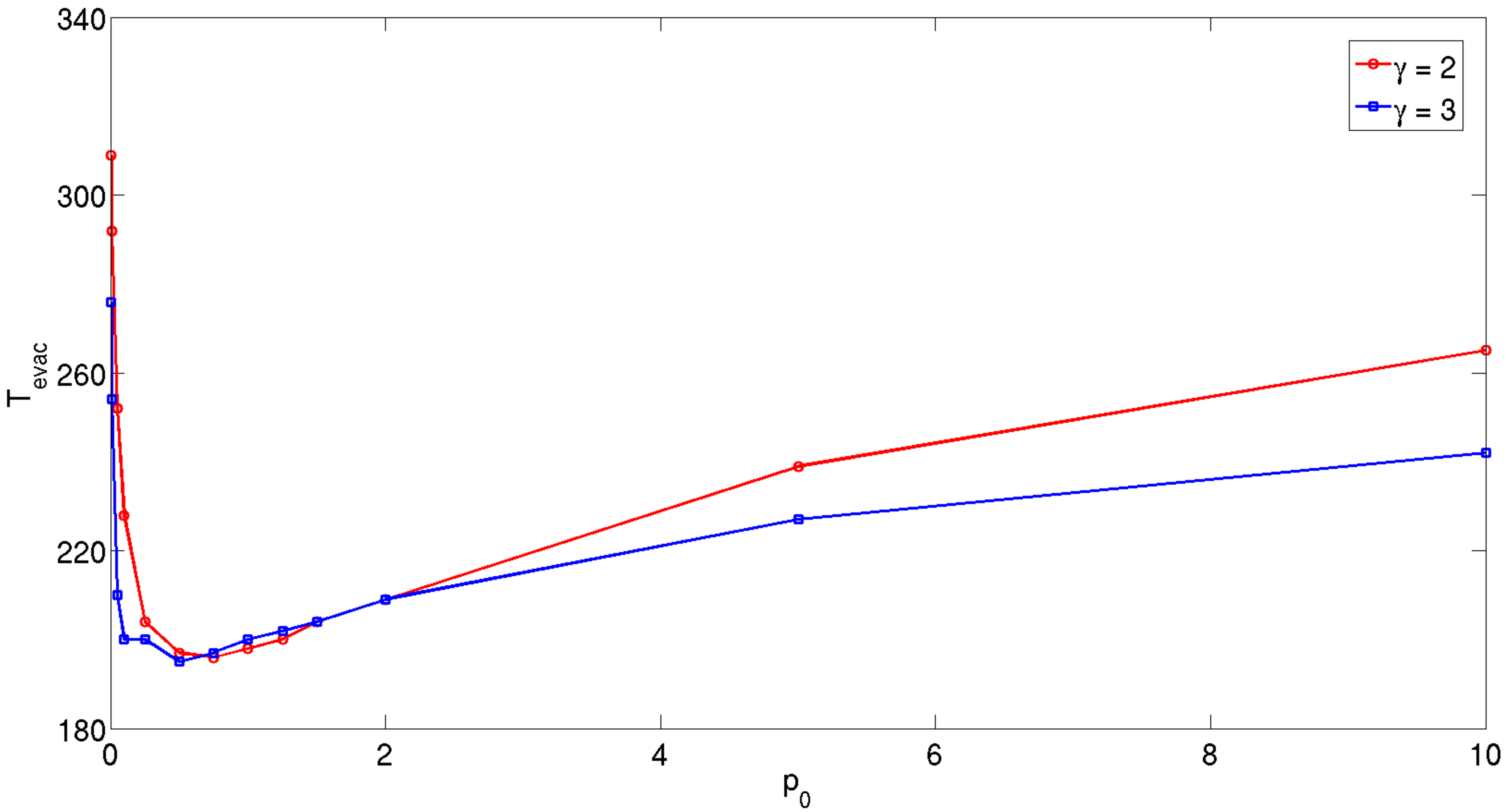}
\end{tabular}
\caption{{\emph{\small{Total evacuation time $T_{\textrm{evac}}$ for the second order model \eqref{eq:mainSystem} with $v_{\max}=2$ m/s, $\rho_{\max}=7$ $\textrm{ped/m}^2$, $\gamma=2,3$ as a function of the pressure coefficient $p_{0}=\{5\times 10^{-3},10^{-2},5\times 10^{-2},10^{-1},2.5\times 10^{-1},5\times 10^{-1},7.5\times 10^{-1},1,1.25,1.5,2,5,10\}$. The initial density is $\rho_{0}=1.5$ $\textrm{ped/m}^2$ in $\Omega_{0}=[1\textrm{ m},5\textrm{ m}]\times[1\textrm{ m},5\textrm{ m}]$, and the initial velocity $\vec{v}_{0}=0$ m/s.} }}}
\label{fig:pressureOptimal}
\end{center}
\end{figure}

 We consider an empty room $10\textrm{ m}\times 6\textrm{ m}$ with a $1$ m wide, symmetrically placed exit and set the following initial data $\rho_{0}=1.5$ $\textrm{ped/m}^2$ in $\Omega_{0}=[1\textrm{ m},5\textrm{ m}]\times[1\textrm{ m},5\textrm{ m}]$, $\vec{v}_{0}=0$ m/s. In Fig~\ref{fig:pressure}a we present the time evolution of the total mass of pedestrians remaining in the room for different values of the coefficient $p_{0}$ in the cases of adiabatic exponent $\gamma=2$ and $v_{\max}=2$ m/s, $\rho_{\max}=7$ $\textrm{ped/m}^2$. We observe that for small values of $p_{0}$ there is a significant decrease of the outflow of pedestrians through the exit, which means that it is blocked. Increasing $p_{0}$ prevents from congestion and as a result the outflow is smother. However, when the distance between pedestrians increases, leaving the room becomes more time consuming. Fig~\ref{fig:pressureOptimal} shows the dependence of the total evacuation time $T_{\textrm{evac}}$ on the parameter $p_{0}$ in the case of $\gamma=2,3$. We observe an optimal value of the parameter $p_{0}\sim 0.5$, which minimizes the evacuation time of pedestrians from the room. 
 
The effect of the adiabatic exponent $\gamma$ is similar. High values increase the repealing forces between pedestrians. We present in Fig~\ref{fig:pressure}b the evolution of the total mass for $\gamma=2,3,4,5$. Clogging observed in the case $\gamma=2$ diminishes for larger $\gamma$. For $\gamma=5$ we observe quasi-linear decrease of the total mass.
 
 The response of pedestrians to compression has an essential effect on the evacuation time. In  normal conditions, when pedestrians want to keep a certain comfort and have enough free space to move, the outflow through the exit is undisturbed. During emergency situations the distances between individuals decrease and the density of the crowd increases, reducing the mobility of pedestrians. As a result the flow may become discontinuous and exits may be blocked. Our simulations show that the second order model is able to reproduce these phenomena. 
 
\subsubsection{Dependence on the desired speed $v_{\max}$}
\begin{figure}[htbp!]
\begin{center}
\begin{tabular}{cc}
\includegraphics[scale=0.11]{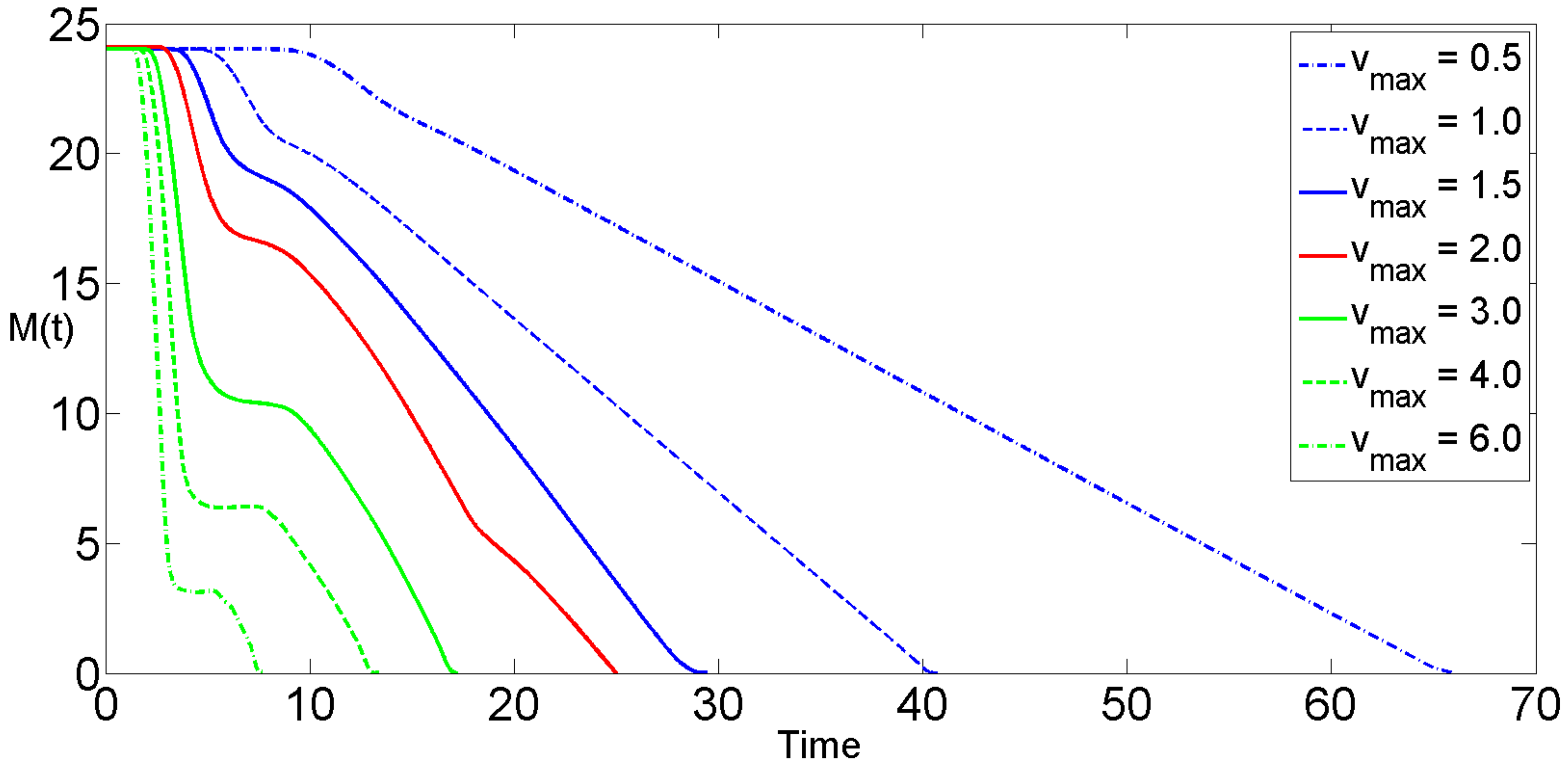}&\includegraphics[scale=0.11]{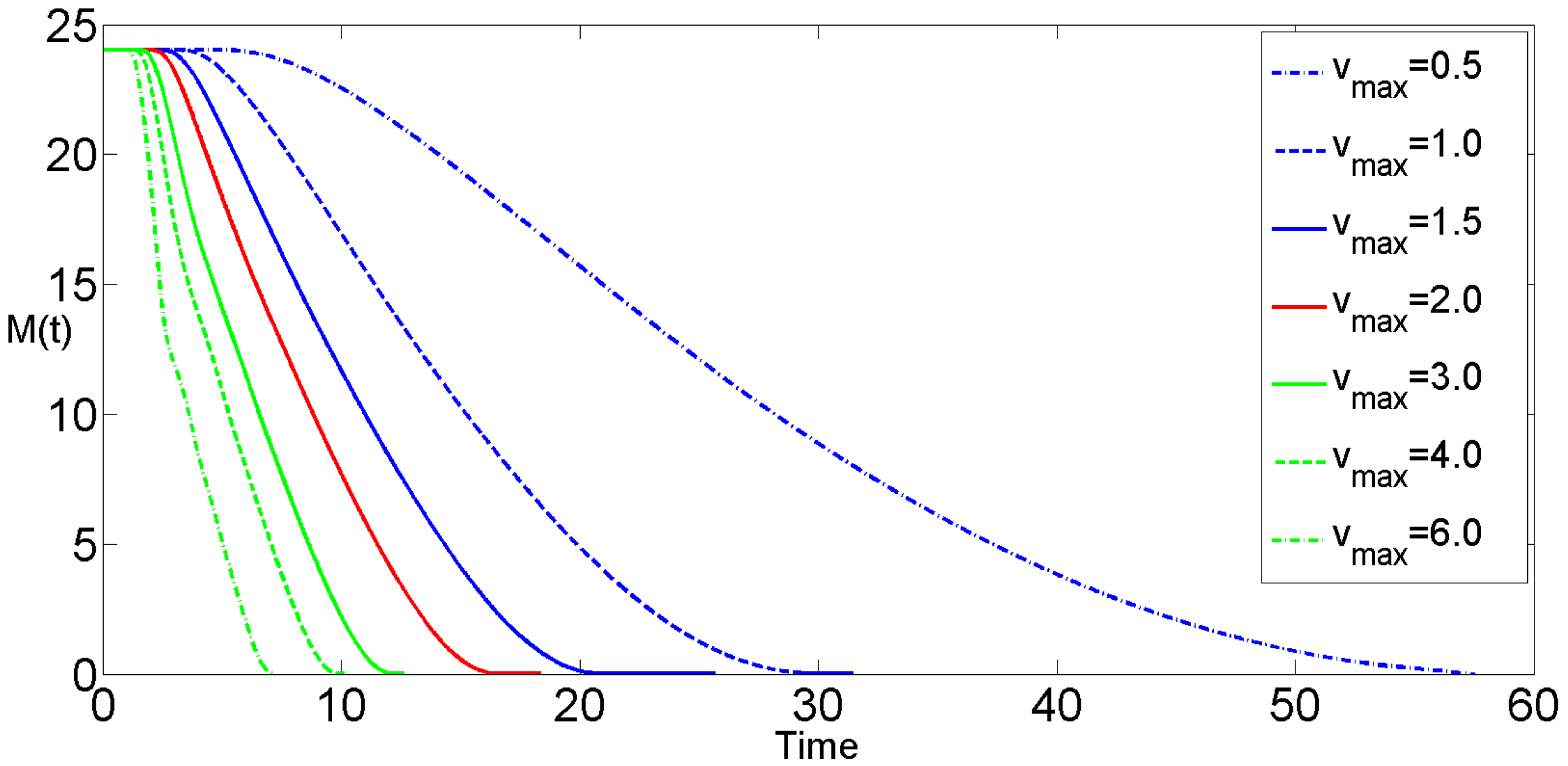}
\end{tabular}
\caption{{\emph{\small{Time evolution of the total mass $M(t)$ of pedestrians inside a room without obstacles for the second order model \eqref{eq:mainSystem} with $\rho_{\max}=7$ $\textrm{ped/m}^2$, $\gamma=2$ and different values of the desired velocity $v_{\max}=\{0.5, 1, 1.5, 2, 3, 4, 6\}$ m/s and different pressure coefficients $p_{0}=0.005$ (on the left) and $p_{0}=0.5$ (on the right). The initial density is $\rho_{0}=5$ $\textrm{ped/m}^2$ in $\Omega_{0}=[1\textrm{ m},5\textrm{ m}]\times[1\textrm{ m},5\textrm{ m}]$, and the initial velocity $\vec{v}_{0}=0$ m/s.} }}}
\label{fig:velocity}
\end{center}
\end{figure}

Now we look for the dependence of the total evacuation time on the desired walking speed $v_{\max}$ of pedestrians. Empirical studies indicate that the average free speed, that is the speed at which pedestrians walk when they are not influenced by others, is about $1.34$ m/s with standard deviation of $0.37$ m/s (see \cite{Buchmueller}). Due to impatience, emergency or panic, people tend to move faster to escape uncomfortable situation or direct life thread as soon as possible. Using the social force model \cite{Helbing95} Helbing et al. \cite{Helbing_Farkas_Vicsek} analyzed an evacuation of 200 people from a room for different desired speeds corresponding to different states of panic. Under the condition of high friction due to the tangential motion of pedestrians, they observed the existence of an optimal speed for which the evacuation was optimized. In Fig~\ref{fig:velocity} we analyze the evacuation of pedestrians from the empty room $\Omega=[0,10\textrm{ m}]\times[0,6\textrm{ m}]$ through a $1$ m  wide exit with $\rho_{0}=2$ $\textrm{ped/m}^2$ in $\Omega_{0}=[1\textrm{ m},5\textrm{ m}]\times[1\textrm{ m},5\textrm{ m}]$, $\vec{v}_{0}=0$ m/s for different values of the desired speed $v_{\max}=0.5,1,1.5,2,3,4,6$ m/s in the case of $p_{0}=0.005$ (on the left) and $p_{0}=0.5$ (on the right).

We observe that unlike Helbing et al. \cite{Helbing_Farkas_Vicsek}, the total evacuation time decreases with higher desired speed. Pedestrians at the front of a group can move almost at their desired velocity $v_{\max}$ as they are not slowed down by the presence of others. When the free speed is high, they reach the exit very fast and leave the room in a short time. However, due to the limited flow capacity of the exit, in the case of $p_{0}=0.005$ we see that pedestrians start to accumulate in front of the door and block it. Decreasing the value of $v_{\max}$ the flow through the exit decreases, so clogging occurs earlier than in the case of larger desired speed. At the same time, also the accumulation has a smaller rate, due to smaller $v_{\max}$,  so the outflow is only slowed down instead of being blocked. Increasing the value of the internal pressure coefficient to $p_{0}=0.5$, the evacuation becomes faster and more regular, as we have already observed in the previous simulations. 

\subsubsection{Effect of obstacles}
In this section we study the evacuation from a room following the idea of Hughes \cite{Hughes2003}, who raised the question of weather suitably placed obstacles can increase the flow through an exit. This idea is an inversion of the Braess paradox \cite{Braess1968, Braess2005}, which was formulated for traffic flows and states that adding extra capacity to a network can in some cases reduce the overall performance. In the case of crowd dynamics, placing an obstacle may be seen intuitively as a worse condition. Nevertheless, it is expected to decrease the density in front of the exit and as a result prevent it from blocking. 

We examine closely this phenomenon in case of pedestrian motion using the second order model \eqref{eq:mainSystem} with $v_{\max}=2$ m/s, $\rho_{\max}=7$ $\textrm{ped/m}^2$. In Section~\ref{sec:effectOfObstacle} we have already observed that an obstacle in front of an exit can reduce the clogging. In this section we consider a room $10\textrm{ m}\times 6\textrm{ m}$ with a $1$ m wide, symmetrically placed exit and system of five circular columns arranged in the shape of a triangle opened towards the door. The columns are centered at $(9.5,2), (9,2.5), (8.5,3), (9,3.5), (9.5,4)$ and have the radius $r=0.22$ m, see Fig~\ref{fig:eikonalColumns} and the initial density equals $\rho_{0}=1$ $\textrm{ped/m}^2$ in $\Omega_{0}=[1\textrm{ m},5\textrm{ m}]\times [1\textrm{ m},5\textrm{ m}]$.

\begin{figure}[htbp!]
\begin{center}
\begin{tabular}{cc}
\includegraphics[scale=0.11]{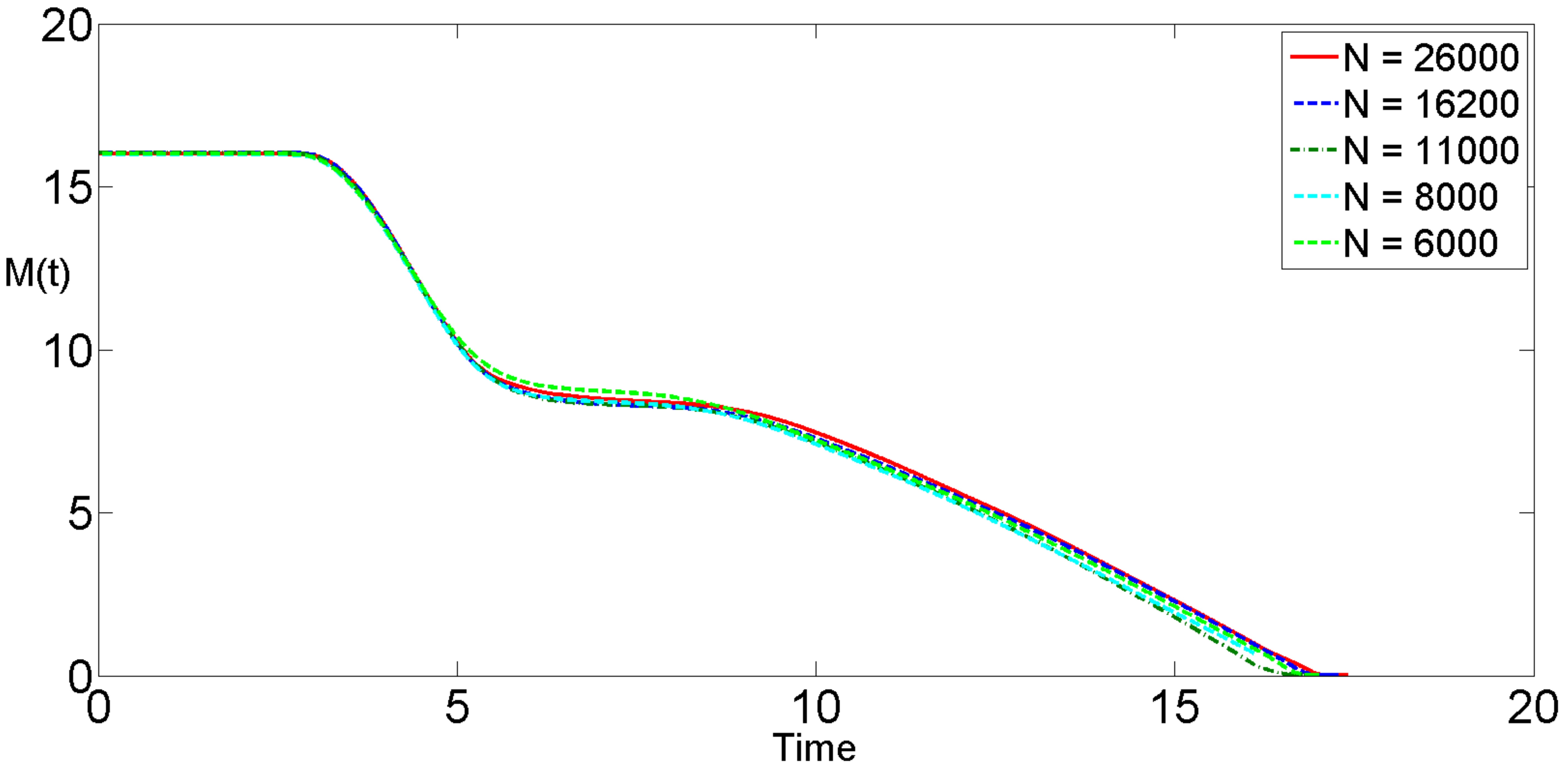}&
\includegraphics[scale=0.11]{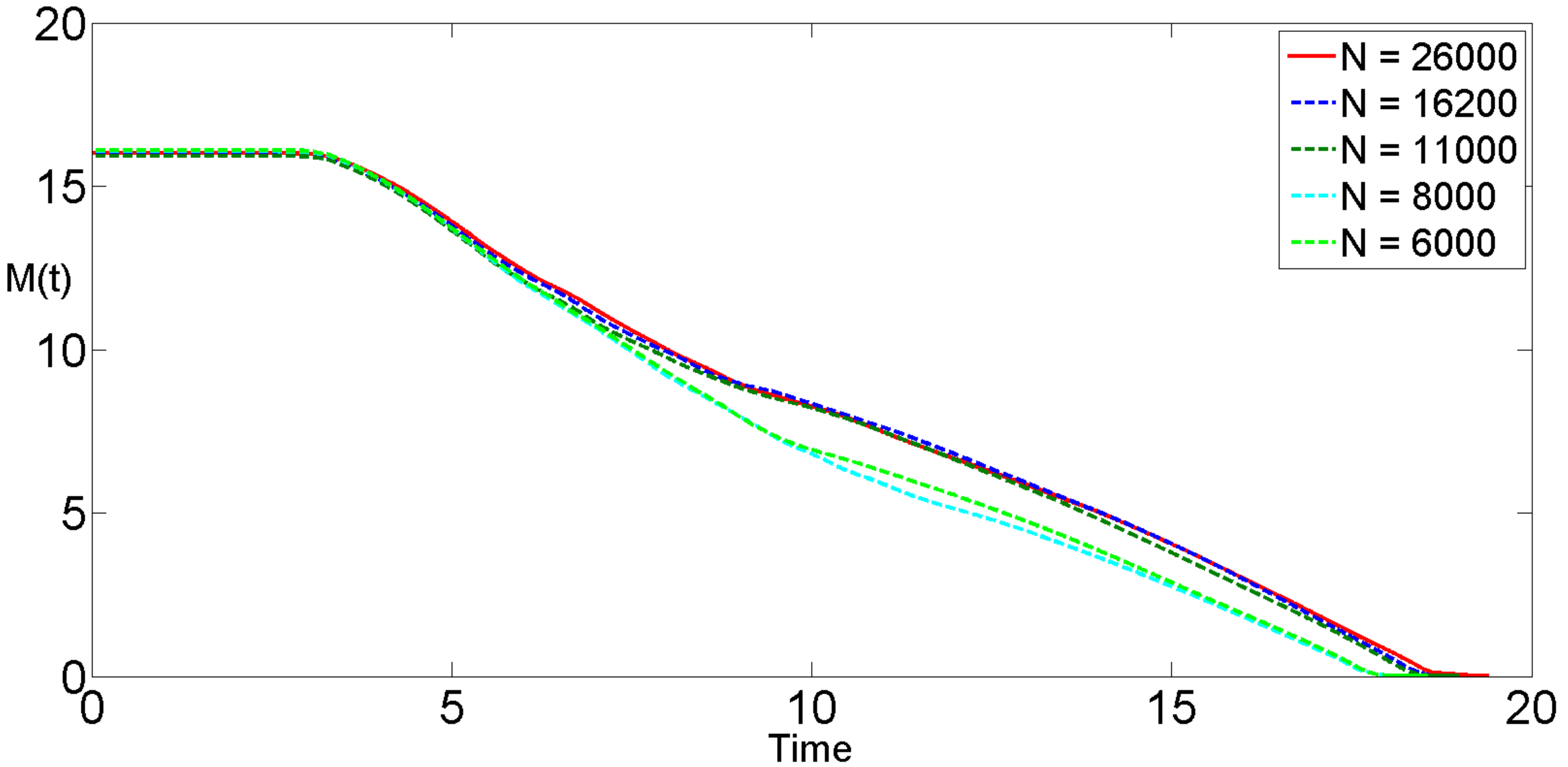}
\end{tabular}
\caption{{\emph{\small{Time evolution of the total mass $M(t)$ of pedestrians for the second order model  \eqref{eq:mainSystem} with $v_{\max}=2$ m/s, $\rho_{\max}=7$ $\textrm{ped/m}^2$, $p_{0}=0.005$, $\gamma=2$ for different number $N$ of finite volume cells in case of an empty room (on the left) and a room with five columns in front of the exit (on the right).} }}}
\label{fig:totalMassMesh}
\end{center}
\end{figure}

\begin{figure}[htbp!]
\begin{center}
\begin{tabular}{cc}
\includegraphics[scale=0.25]{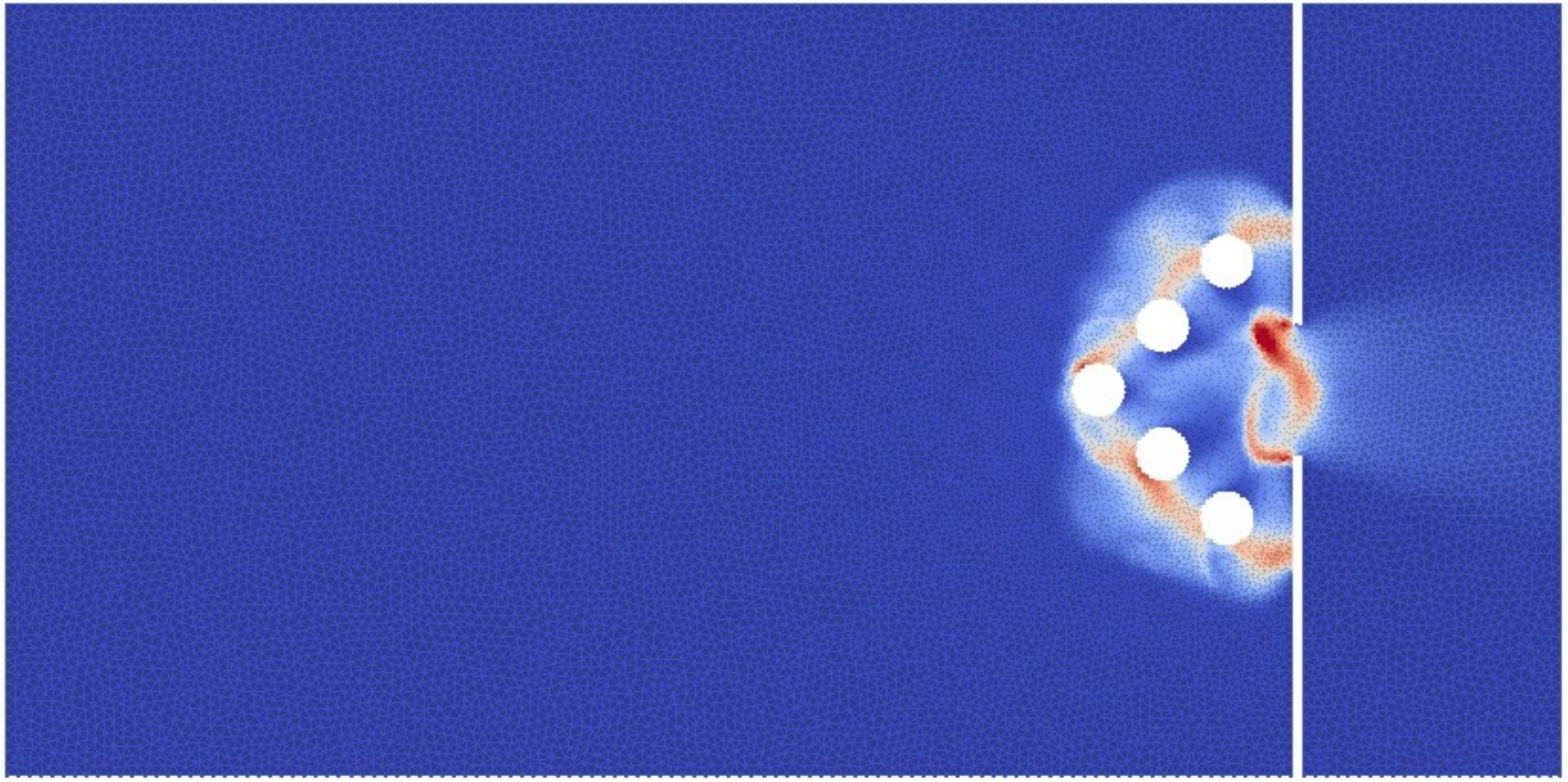}&
\includegraphics[scale=0.25]{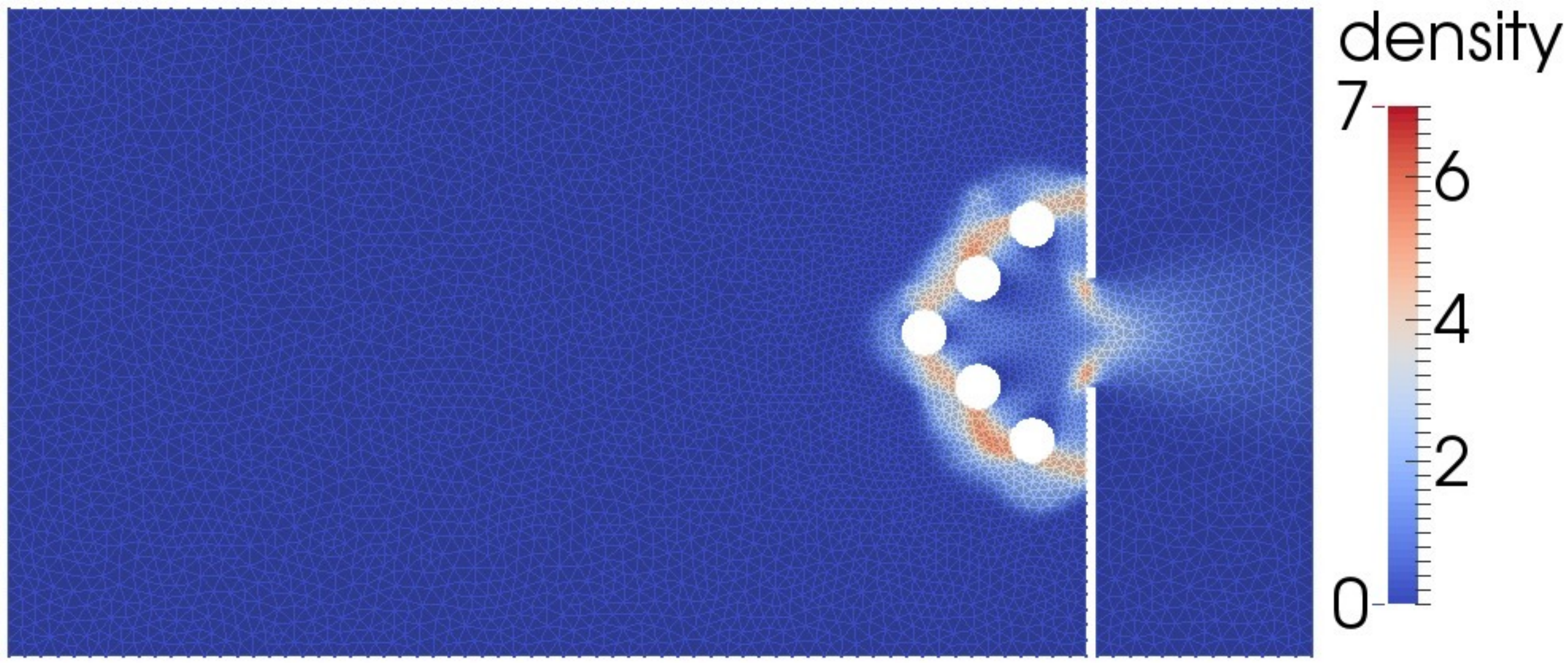}\\
fine mesh&coarse mesh
\end{tabular}
\caption{{\emph{\small{Density profiles at time $T=9$ s for the second order model \eqref{eq:mainSystem} with $v_{\max}=2$ m/s, $\rho_{\max}=7$ $\textrm{ped/m}^2$, $p_{0}=0.005$, $\gamma=2$ for different numbers of the finite volume cells: $N = 16000$ (on the left) and $N = 6000$ (on the right). }}}}
\label{fig:lengthScale}
\end{center}
\end{figure}

In order to compare the efficiency of the evacuation we use the total mass of pedestrians that remain in the room, which corresponds directly to the outflow through the door. Fig~\ref{fig:totalMassMesh} presents the time evolution of the total mass on meshes with different number of finite volume cells, which are refined by a factor two near the exit and the obstacles. The room with the columns and the obstacle-free case are considered. We observe that for the empty room (on the left) the results do not change with the mesh refinement. On the contrary, in the case with the obstacles (on the right) there is a significant difference between the total mass curves for $N=16000,26000$ and for $N=6000,8000$ finite volume cells. This effect is a result of the presence of phenomena occurring at different length scales for irregular flows. In Fig~\ref{fig:lengthScale} we can see that the clogging observed on the mesh with $N=16000$ finite volume cells is not resolved for $N=6000$ resulting in false faster evacuation.  

\begin{figure}[htbp!]
\begin{center}
\begin{tabular}{cc}
\includegraphics[scale=0.11]{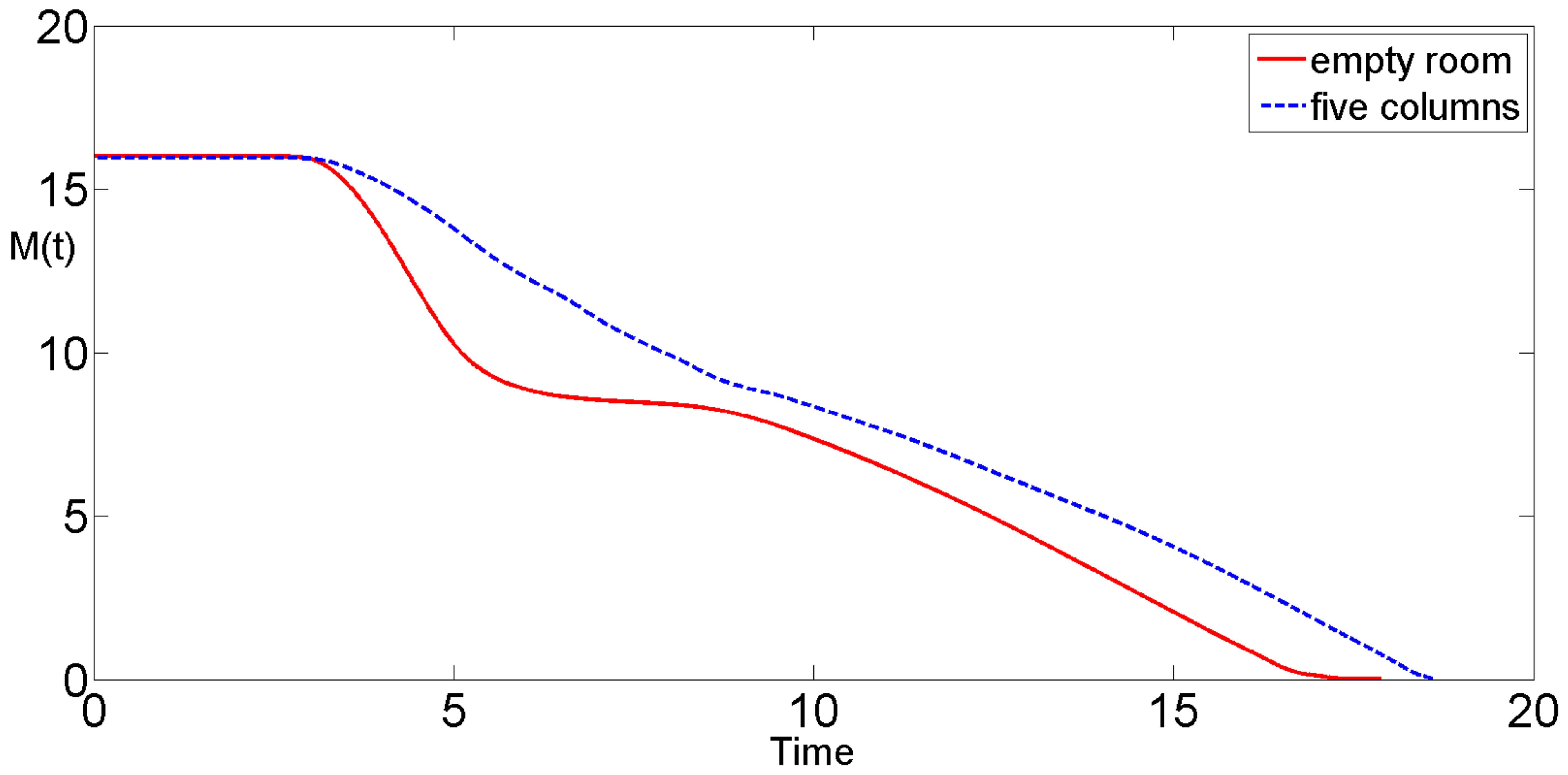}&
\includegraphics[scale=0.11]{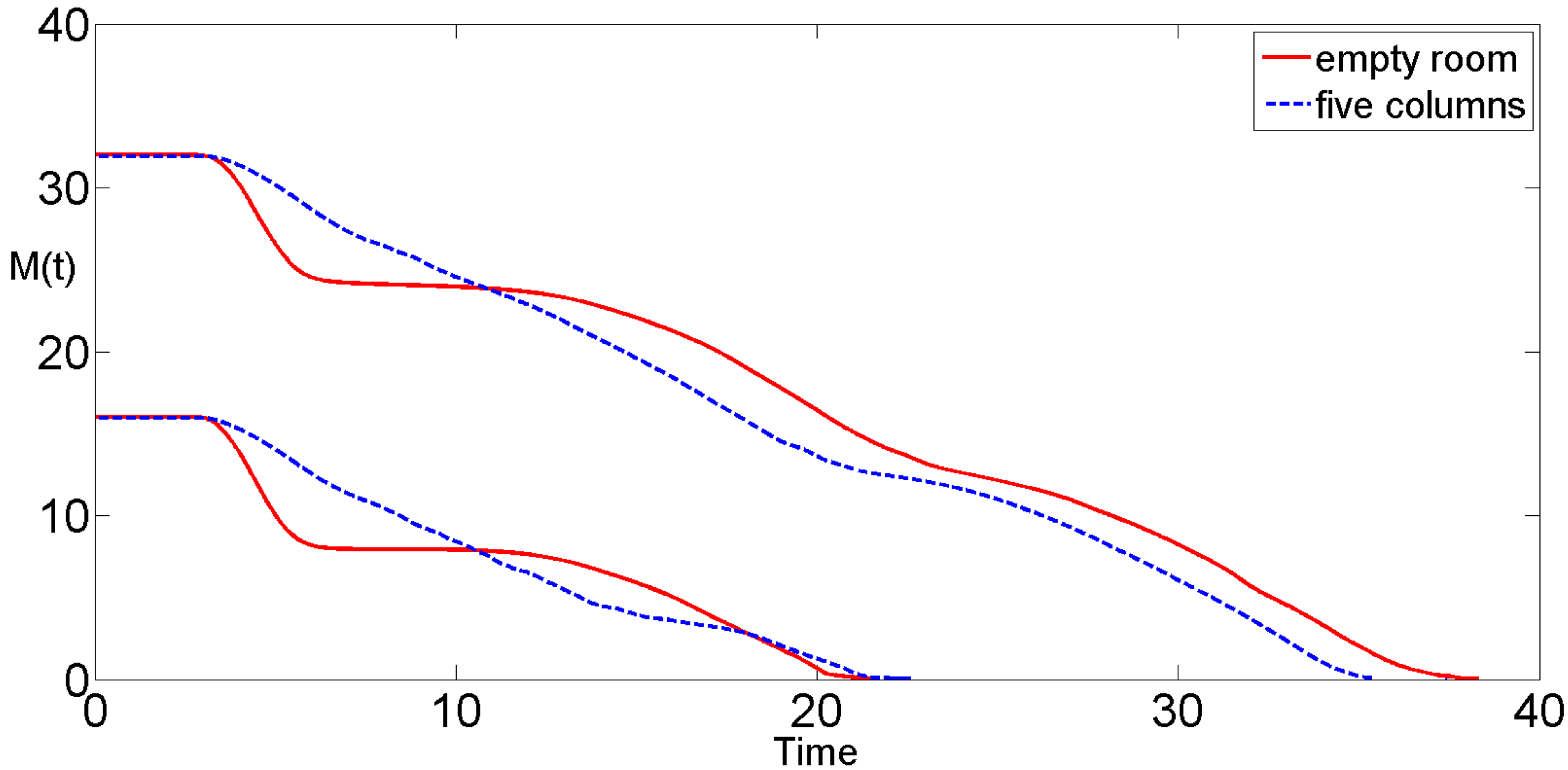}
\end{tabular}
\caption{{\emph{\small{Time evolution of the total mass $M(t)$ of pedestrians for the second order model \eqref{eq:mainSystem} with $v_{\max}=2$ m/s, $\rho_{\max}=7$ $\textrm{ped/m}^2$, $\gamma=2$ and different pressure coefficients: $p_{0}=0.005$ (on the left) and $p_{0}=0.001$ (on the right). }}}}
\label{fig:evacuationObstacles}
\end{center}
\end{figure}

In the following test we compare numerically the evacuation from the room described above using a mesh with $N=16000$ finite volume cells. Fig~\ref{fig:evacuationObstacles} shows the time evolution of the total mass for two different pressure coefficients $p_{0}=0.005,0.001$. We observe that the clogging present in the empty room is reduced significantly using the obstacles. The system of columns creates a ''waiting zone'' in front of the exit. Pedestrians are slowed down and partially stopped by the obstacles, which corresponds to the slower outflow in the initial phase of the evacuation. But, at the same time the density at the exit remains low because the incoming flow does not exceed the door capacity. In case of small pressure $p_{0}=0.001$, which allows for higher congestion, the improvement due to the obstacles is more visible. For the initial mass $M=32$ and columns of radius $r=0.24$ m the total evacuation time of the room with the columns becomes smaller that in the obstacles-free situation.   

\section{Conclusions}
In this study, two pedestrian flow models have been analyzed in the context of macroscopic modelling of evacuation of pedestrians from a room through a narrow exit. Error analysis has been carried out for numerical validation of a finite-volume scheme on unstructured grid.Various test-cases have been considered to quantify the influence of the model parameters on the behaviour of solutions and to measure the ability of the models to reproduce some of the phenomena occurring in the evacuation of high density crowds. More precisely, numerical experiments show that the classical Hughes' type model cannot reproduce stop-and-go waves or clogging at bottlenecks. On the other hand, it was verified numerically that the second order model captures better the structure of interactions between pedestrians and is able to produce the above behaviours. However, even this model is still far from being validated and should be verified and calibrated with realistic experiments. In fact, we have pointed out that  values of some of its parameters have a significant effect on the formation of the above phenomena so their tuning is essential. For some particular choices of the parameters the evacuation through a narrow exit was analyzed and an example of the inverse Braess paradox was given. It was shown that using a particular configuration of obstacles it is possible to reduce the clogging at the exit and increase the outflow. Analysis of more realistic settings are to be considered at the next step. 

\section*{Acknowledgment}
This research was supported by  the European Research Council
under the European Union's Seventh Framework Program (FP/2007-2013) / ERC Grant Agreement n. 257661.

\bibliographystyle{plain}

\begin{thebibliography}{10}

\bibitem{Bellomo_Bellouquid}
N.~Bellomo and A.~Bellouquid.
\newblock On the modelling of vehicular traffic and crowds by kinetic theory of
  active particles.
\newblock In {\em Mathematical Modeling of Collective Behavior in
  Socio-Economic and Life Sciences}, number~2 in Modeling and Simulation in
  Science, Engineering and Technology, pages 215--221. Birkhäuser Boston,
  2010.

\bibitem{Bellomo}
N.~Bellomo and C.~Dogb\'{e}.
\newblock On the modelling crowd dynamics from scaling to hyperbolic
  macroscopic models.
\newblock {\em Mathematical Models and Methods in Applied Sciences},
  18:1317--1345, 2008.

\bibitem{Bornemann_Rasch}
F.~Bornemann and Ch. Rasch.
\newblock Finite-element discretization of static {H}amilton-{J}acobi equations
  based on a local variational principle.
\newblock {\em Computing and Visualization in Science}, 9(2):57--69, 2006.

\bibitem{Braess1968}
D.~Braess.
\newblock {\"{U}ber} ein {Paradoxon} aus der {Verkehrsplanung}.
\newblock {\em Unternehmensforschung}, 12:258--268, 1968.

\bibitem{Braess2005}
D.~Braess, A.~Nagurney, and T.~Wakolbinger.
\newblock On a paradox of traffic planning.
\newblock {\em Transportation Science}, 39(4):446--450, 2005.

\bibitem{Bryson_Levy}
S.~Bryson and D.~Levy.
\newblock High-order semi-discrete central-upwind schemes for multi-dimensional
  {H}amilton-{J}acobi equations.
\newblock {\em Journal of Computational Physics}, 189(1):63--87, 2003.

\bibitem{Buchmueller}
S.~Buchmueller and U.~Weidmann.
\newblock {\em Parameters of Pedestrians, Pedestrian Traffic and Walking
  Facilities}.
\newblock Schriftenreihe des IVT. ETH Zurich, 2006.

\bibitem{ColomboRosini}
R.M. Colombo and M.D. Rosini.
\newblock Pedestrian flows and nonclassical shocks.
\newblock {\em Mathematical Methods in the Applied Sciences},
  28(13):1553--1567, 2008.

\bibitem{CristianiPiccoliTosin}
Emiliano Cristiani, Benedetto Piccoli, and Andrea Tosin.
\newblock Multiscale modeling of granular flows with application to crowd
  dynamics.
\newblock {\em Multiscale Modeling \& Simulation}, 9(1):155--182, 2011.

\bibitem{Degond_etal}
P.~Degond, C.~Appert-Rolland, J.~Pettr\'{e}, and G.~Theraulaz.
\newblock Vision-based macroscopic pedestrian models.
\newblock {\em Kinetic and Related Models}, 6(4):809--839, 2013.

\bibitem{Dijkstra}
J.~Dijkstra, H.~J.~P. Timmermans, and A.~J. Jessurun.
\newblock A multi-agent cellular automata system for visualising simulated
  pedestrian activity.
\newblock pages 29--36. Springer Verlag, 2000.

\bibitem{Einfeldt}
B.~Einfeldt, C.~D. Munz, P.~L. Roe, and B.~Sj\"{o}green.
\newblock On {G}odunov-type methods near low densities.
\newblock {\em Journal of Computational Physics}, 92(2):273--295, 1991.

\bibitem{Escobar}
R.~Escobar and A.~De~La~Rosa.
\newblock Architectural design for the survival optimization of panicking
  fleeing victims.
\newblock volume 2801, pages 97--106. Springer, 2003.

\bibitem{Falcone_Ferretti}
M.~Falcone and R.~Ferretti.
\newblock Semi-{L}agrangian schemes for {H}amilton-{J}acobi equations, discrete
  representation formulae and {G}odunov methods.
\newblock {\em Journal of Computational Physics}, 175(2):559--575, 2002.

\bibitem{Frank_Dorso}
G.A. Frank and C.O. Dorso.
\newblock Room evacuation in the presence of an obstacle.
\newblock {\em Physica A: Statistical Mechanics and its Applications},
  390(11):2135--2145, 2011.

\bibitem{Gopal}
S.~Gopal and T.~R. Smith.
\newblock Navigator: an ai-based model of human way-finding in an urban
  environment.
\newblock {\em Spatial Choices and Processes}, pages 169--200, 1989.

\bibitem{HLL}
A.~Harten, P.~D. Lax, and B.~van Leer.
\newblock On upstream differencing and godunov-type schemes for hyperbolic
  conservation laws.
\newblock {\em SIAM Review}, 25:35--61, 1983.

\bibitem{Helbing_Buzna2005}
D.~Helbing, L.~Buzna, A.~Johansson, and Torsten Werner.
\newblock Self-organized pedestrian crowd dynamics: Experiments, simulations,
  and design solutions.
\newblock {\em Transportation Science}, 39:1--24, 2005.

\bibitem{Helbing_Farkas_Vicsek}
D.~Helbing, I.~Farkas, P.~Moln\`{a}r, and T.~Vicsek.
\newblock Simulation of pedestrian crowds in normal and evacuation situations.
\newblock In {\em Pedestrian and Evacuation Dynamics}, pages 21--58. Springer,
  2002.

\bibitem{Helbing_Nature}
D.~Helbing, I.~Farkas, and T.~Vicsek.
\newblock Simulating dynamical features of escape panic.
\newblock {\em Nature}, 407:487--490, 2000.

\bibitem{Helbing2002}
D.~Helbing, I.~Farkas, and T.~Vicsek.
\newblock Crowd disasters and simulation of panic situations.
\newblock In {\em The Science of Disasters}, pages 330--350. Springer Berlin
  Heidelberg, 2002.

\bibitem{Helbing_freezing}
D.~Helbing, I.J. Farkas, and T.~Vicsek.
\newblock Freezing by heating in a driven mesoscopic system.
\newblock {\em Physical Review Letters}, 84:1240--1243, 2000.

\bibitem{Helbing_Johansson2009}
D.~Helbing and A.~Johansson.
\newblock Pedestrian, crowd and evacuation dynamics.
\newblock In {\em Encyclopedia of Complexity and Systems Science}, pages
  6476--6495. Springer, 2009.

\bibitem{Helbing_Johansson_Abideen}
D.~Helbing, A.~Johansson, and H.~Zein Al-Abideen.
\newblock Dynamics of crowd disasters: An empirical study.
\newblock {\em Physical Review E}, 75(4):046109, 2007.

\bibitem{Helbing95}
D.~Helbing and P.~Moln\`{a}r.
\newblock Social force model for pedestrian dynamics.
\newblock {\em Physical Review E}, pages 4282--4286, 1995.

\bibitem{HoogendoornDaamen}
S.P. Hoogendoorn and W.~Daamen.
\newblock Pedestrian behavior at bottlenecks.
\newblock {\em Transportation Science}, 39(2):147--159, 2005.

\bibitem{Hu_Shu}
Ch. Hu and Ch.W. Shu.
\newblock A discontinuous {G}alerkin finite element method for
  {H}amilton-{J}acobi equations.
\newblock {\em SIAM Journal on Scientific Computing}, 21(2):666--690, 1999.

\bibitem{Huang}
L.~Huang, S.C. Wong, M.~Zhang, C.-W. Shu, and W.H.K. Lam.
\newblock Revisiting hughes' dynamic continuum model for pedestrian flow and
  the development of an efficient solution algorithm.
\newblock {\em Transportation Research Part B: Methodological}, 43(1):127--141,
  2009.

\bibitem{Hughes2002}
R.~L. Hughes.
\newblock A continuum theory for the flow of pedestrians.
\newblock {\em Transportation Research Part B: Methodological}, 36(6):507 --
  535, 2002.

\bibitem{Hughes2003}
R.L. Hughes.
\newblock The flow of human crowds.
\newblock In {\em Annual review of fluid mechanics}, volume~35, pages 169--182.
  2003.

\bibitem{num3sis}
Inria.
\newblock Num3sis software, 2012.

\bibitem{Jiang_Liu}
Y.-Q Jiang, R.-X. Liu, and Y.-L. Duan.
\newblock Numerical simulation of pedestrian flow past a circular obstruction.
\newblock {\em Acta Mechanica Sinica}, 27(2):215--221, 2011.

\bibitem{Jiang2010}
Y.Q. Jiang, P.~Zhang, S.C. Wong, and R.X. Liu.
\newblock A higher-order macroscopic model for pedestrian flows.
\newblock {\em Physica A: Statistical Mechanics and its Applications},
  389(21):4623 -- 4635, 2010.

\bibitem{Keating}
J.~P. Keating.
\newblock The myth of panic.
\newblock {\em Fire Journal}, May:57--62, 1982.

\bibitem{Kretz}
T.~Kretz, A.~Gr\"{u}nebohm, and M.~Schreckenberg.
\newblock Experimental study of pedestrian flow through a bottleneck.
\newblock {\em Journal of Statistical Mechanics: Theory and Experiment},
  (10):10014, 2006.

\bibitem{num3sisOpale}
J.~Labroqu\`{e}re, R.~Duvigneau, T.~Kliczko, and J.~Wintz.
\newblock Interactive computation and visualization towards a virtual wind
  tunnel.
\newblock {\em 47th 3AF Symposium on Applied Aerodynamics, Paris, France,},
  March, 2012.

\bibitem{LachapelleWolfram}
Aimé Lachapelle and Marie-Therese Wolfram.
\newblock On a mean field game approach modeling congestion and aversion in
  pedestrian crowds.
\newblock {\em Transportation Research Part B: Methodological}, 45(10):1572 --
  1589, 2011.

\bibitem{Lemercier}
S.~Lemercier, A.~Jelic, R.~Kulpa, J.~Hua, J.~Fehrenbach, P.~Degond,
  C.~Appert-Rolland, S.~Donikian, and J.~Pettr\'{e}.
\newblock Realistic following behaviors for crowd simulation.
\newblock {\em Computer Graphics Forum}, 31:489--498, 2012.

\bibitem{Maury}
B.~Maury, A.~Roudneff-Chupin, and F.~Santambrogio.
\newblock A macroscopic crowd motion model of gradient flow type.
\newblock {\em Mathematical Models and Methods in Applied Sciences},
  20:1787--1821, 2009.

\bibitem{Maury_etal}
Bertrand Maury, Aude Roudneff-Chupin, Filippo Santambrogio, and Juliette Venel.
\newblock Handling congestion in crowd motion modeling.
\newblock {\em Networks and Heterogeneous Media}, 6(3):485--519, 2011.

\bibitem{Muramatsu}
M.~Muramatsu, T.~Irie, and T.~Nagatani.
\newblock Jamming transition in pedestrian counter flow.
\newblock {\em Physica A: Statistical and Theoretical Physics},
  267(3-4):487--498, 1999.

\bibitem{Osher_Fedkiw}
S.~Osher and R.~Fedkiw.
\newblock Level set methods and dynamic implicit surfaces.
\newblock {\em Applied Mathematical Sciences}, 153, 2003.

\bibitem{Osher_Sethian}
S.~Osher and J.~A. Sethian.
\newblock Fronts propagating with curvature-dependent speed: algorithms based
  on {H}amilton-{J}acobi formulations.
\newblock {\em Journal of Computational Physics}, 79(1):12--49, 1988.

\bibitem{Payne1971}
H.J. Payne.
\newblock {\em Models of Freeway Traffic and Control}.
\newblock Simulation Councils, Incorporated, 1971.

\bibitem{Tosin}
B.~Piccoli and A.~Tosin.
\newblock Pedestrian flows in bounded domains with obstacles.
\newblock {\em Continuum Mechanics and Thermodynamics}, 21(2):85--107, 2009.

\bibitem{predtechenskii1978planning}
V.M. Predtechenskii and A.I. Milinskii.
\newblock {\em Planning for Foot Traffic Flow in Buildings}.
\newblock TT. Amerind, 1978.

\bibitem{Roache}
P.~J. Roache.
\newblock Perspective: A method for uniform reporting of grid refinement
  studies.
\newblock {\em Journal of Fluids Engineering}, 116(3):405--413, 1994.

\bibitem{Sethian}
J.~A. Sethian.
\newblock Level set methods and fast marching methods.
\newblock {\em Cambridge Monographs on Applied and Computational Mathematics},
  3, 1999.

\bibitem{Seyfried}
A.~Seyfried, M.~Boltes, J.~Kähler, W.~Klingsch, A.~Portz, T.~Rupprecht,
  A.~Schadschneider, B.~Steffen, and A.~Winkens.
\newblock Enhanced empirical data for the fundamental diagram and the flow
  through bottlenecks.
\newblock In {\em Pedestrian and Evacuation Dynamics}, pages 145 -- 156.
  Berlin/Heidelberg, Springer, 2010.

\bibitem{Seyfried_bottleneck}
A.~Seyfried, O.~Passon, B.~Steffen, M.~Boltes, T.~Rupprecht, and W.~Klingsch.
\newblock New insights into pedestrian flow through bottlenecks.
\newblock {\em Transportation Science}, 43(3):395--406, 2009.

\bibitem{Seyfried2006}
A.~Seyfried, B.~Steffen, and T.~Lippert.
\newblock Basics of modelling the pedestrian flow.
\newblock {\em Physica A}, 368:232--238, 2006.

\bibitem{toro}
E.F. Toro.
\newblock {\em Riemann Solvers and Numerical Methods for Fluid Dynamics: A
  Practical Introduction}.
\newblock Springer-Verlag Berlin Heidelberg, 2009.

\bibitem{Tsai_FSM}
Y.~R. Tsai, L.-T. Cheng, S.~Osher, and H.K. Zhao.
\newblock Fast sweeping algorithms for a class of {H}amilton-{J}acobi
  equations.
\newblock {\em SIAM Journal on Numerical Analysis}, 41(2):673--694, 2003.

\bibitem{InriaReport}
Monika Twarogowska, Paola Goatin, and R{\'e}gis Duvigneau.
\newblock {Numerical study of macroscopic pedestrian flow models}.
\newblock INRIA Research Report no. 8340, July 2013.

\bibitem{Whitham1974}
G.B. Whitham.
\newblock {\em Linear and nonlinear waves}.
\newblock Pure and applied mathematics. Wiley, 1974.

\bibitem{Zuriguel}
I.~Zuriguel, A.~Janda, A.~Garcimart'in, C.~Lozano, R.~Ar\'{e}valo, and D.~Maza.
\newblock Silo clogging reduction by the presence of an obstacle.
\newblock {\em Physical Review Letters}, 107:278001, 2011.

\end{thebibliography}

\end{document}